\documentclass[12pt]{article}

\usepackage[colorlinks=true,
linkcolor=webgreen,
filecolor=webbrown,
citecolor=webgreen]{hyperref}
\definecolor{webgreen}{rgb}{0,.5,0}
\definecolor{webbrown}{rgb}{.6,0,0}

\usepackage{amsmath}
\usepackage{graphics}
\usepackage{amscd}
\usepackage{amssymb}
\def\Reals{\mathop{\hbox{\mit I\kern-.2em R}}\nolimits}

\def\RA{\mathop{\hbox{\Rightarrow\kern-.2em /}}}

\def\Complexes{\mathop{\hbox{\mit C\kern-.44em
               \vrule depth 0ex height 1.4ex width .06em
               \kern.41em}}\nolimits}

\def\Zeals{\mathop{\hbox{\mit Z\kern-.29em Z}}\nolimits}

\def\Neals{\mathop{\hbox{\mit I\kern-.2em N}}\nolimits}

\def\K{\mathop{\hbox{\mit I\kern-.2em K}}\nolimits}

\def\F{\mathop{\hbox{\mit I\kern-.2em F}}\nolimits}

\def\sqr#1#2{{\vcenter{\vbox{\hrule height.#2pt
  \hbox{\vrule width.#2pt height#1pt \kern#1pt
    \vrule width.#2pt}
   \hrule height.2pt}}}}

\def \oer{\overrightarrow{e}}
\def \oel{\overleftarrow{e}}
\def \dbar{\bar{d}}

\def \bG{{\bf G}}

\def \bu{{\bf u}}

\def \bv{{\bf v}}

\def \bX{{\bf X}}

\def \calE{{\cal E}}

\def \calH{{\cal H}}

\newcommand{\la}{\lambda}
\newcommand{\La}{\Lambda}
\newcommand{\af}{\alpha}
\newcommand{\ZZ}{\mathbb Z}

\usepackage[dvips]{epsfig}
\begin{document}
\begin{center}
{\Large \bf Multiple Description Vector Quantization with Lattice Codebooks:
Design and Analysis} \\
\ \\
{\em Vinay A. Vaishampayan, N. J. A. Sloane}\\ 
{\em AT\&T Shannon Laboratory}\\ {\em Florham Park, NJ 07932} \\
\ \\
and\\
\ \\
{\em Sergio D. Servetto}\\
{\em Ecole Polytechnique F\'{e}d\'{e}rale
de Lausanne, CH-1015 Lausanne, Switzerland}

\vspace{0.5\baselineskip}
Jan. 25, 2000; revised Oct. 7, 2000\\
\vspace{0.5\baselineskip}
\end{center}

\begin{center}
{\bf Abstract}
\end{center}
The problem of designing a multiple description vector quantizer
with lattice codebook $\Lambda$ is considered.
A general solution is given to a labeling problem which
plays a crucial role in the design of such quantizers.
Numerical performance results are obtained
for quantizers based on the lattices $A_2$ and $\ZZ^i$, $i=1,2,4,8$,
that make use of this labeling algorithm.

The high-rate squared-error distortions for this family of $L$-dimensional vector
quantizers
are then analyzed for a memoryless source with
probability density function $p$ and differential entropy $h(p) <
\infty$.
For any $a \in (0,1)$ and rate pair $(R,R)$, it is shown that the
two-channel distortion $\bar{d}_0$
and the channel 1 (or channel 2) distortion $\bar{d}_s$
satisfy 
\[
\lim_{R\rightarrow \infty} \bar{d}_02^{2R(1+a)}=
\frac{1}{4}G(\Lambda)2^{2h(p)}
\]
and 
\[
\lim_{R\rightarrow \infty} \bar{d}_s 2^{2R(1-a)} = G(S_L)
2^{2h(p)},
\]
where $G(\Lambda)$ is the normalized second moment of a Voronoi
cell of the lattice $\Lambda$ and
$G(S_L)$ is the normalized second moment of a sphere in
$L$ dimensions.

\vspace{\fill}
\vspace*{\fill}
\clearpage

\baselineskip=21pt
\pagestyle{headings}
\section{Introduction}
\label{introduction}
We consider the problem of designing a multiple description vector quantizer for a memoryless
source with 
probability density function $p$, differential entropy $h(p)< \infty$ and the squared-error
distortion measure.  A multiple description vector quantizer encodes vectors for transmission
over a two-channel
communication system. The objective is to send information about the source vector over each
channel
in such a way that good performance is obtained when both channels work and the degradation
is
small when either channel fails. It is assumed that the encoder has no knowledge about the
state
of a channel, i.e., it does not know whether a channel has failed or is working.

The recent interest in the multiple description problem is largely because of the
application to image, video and voice communications over packet networks with nonzero
probability of packet loss. The loss of a packet could be significant
if it results in a loss of a large block of source samples, such as a large image block or
a large block of speech. One way to improve performance is to place different encoded versions of
a given block of source samples into several packets in such a way that if some of these packets are
received, a degraded version of the source block may be recovered. This leads to the
formulation of the multiple description source coding problem.

For the single-description problem, one of the benefits of vector quantization over scalar
quantization
is a reduction in granular distortion. This is because in higher dimensions it is possible to
construct Voronoi cells
that are more ``spherical'' than the hypercube. To be more
specific, uniform scalar quantization coupled with entropy coding is known to have
mean squared error (mse) $\bar{d}(R)$ at entropy $R$ bits/sample 
satisfying~\cite{GPI}
\begin{equation}
\lim_{R\rightarrow \infty}\bar{d}(R)2^{2R}= \frac{2^{2h(p)}}{12};
\label{eqn11}
\end{equation}
whereas if an $L$-dimensional lattice $\Lambda$ is used as a codebook, the distortion
satisfies
\begin{equation}
\lim_{R\rightarrow \infty} \bar{d}(R)2^{2R}=G(\Lambda)2^{2h(p)},
\label{eqn-1.2}
\end{equation}
where $G(\Lambda)$ is the normalized second moment of a Voronoi cell of the lattice.
In dimensions greater than one, lattices exist for which
$G(\Lambda)$ is strictly smaller than
$1/12$.  For example, in 8 dimensions, it is possible to
gain 0.66 dB by using the lattice  $E_8$ as compared to uniform scalar
quantization~\cite{SPLAG}.
It is also known through a random quantizing argument \cite{Zador1} that quantizers exist
for which the product $\bar{d}(R)2^{2R}$ approaches
$2^{2h(p)}/(2\pi e)$
as the rate increases.
Furthermore, it follows from rate distortion theory~\cite{BERI} that no smaller value can be
achieved for the above product in the limit of infinite rate.
The maximum gain possible over entropy-coded scalar
quantization is 1.53 dB and
lattices provide a useful method for closing this  gap.

From now on we will restrict our attention to the case
of two channels.
Consider a multiple description quantizer
which sends information across each channel at 
a rate of $R$ bits/sample.
The performance of this  system  is measured in terms of three distortions:
the two-channel distortion $\bar{d}_0$, when both
descriptions are available to the decoder; the  channel 1 distortion $\bar{d}_1$, 
when only the first description is available and the  channel 2 distortion $\bar{d}_2$, 
when only the second description is available. We will further assume that
$\bar{d}_1=\bar{d}_2=\bar{d}_s$ and will refer to this common value as the side distortion.
The 
objective is to design vector quantizers that minimize $\bar{d}_0$ under the constraint
$\bar{d}_s \leq
D_s$, for a given rate pair $(R,R)$ and a given bound $D_s$ on the side-channel
distortion.

It has been shown \cite{VBII} that for a uniform entropy-coded multiple description
quantizer, and any $a \in (0,1)$, the distortions satisfy
\begin{eqnarray}
\lim_{R\rightarrow \infty}\bar{d}_0(R)2^{2R(1+a)} & = &
\frac{1}{4}\left(\frac{2^{2h(p)}}{12}\right), \nonumber \\
\lim_{R\rightarrow \infty}\bar{d}_s(R)2^{2R(1-a)} & = & \left(\frac{2^{2h(p)}}{12} \right).     
\label{eqn-ecmdsq0}
\end{eqnarray}
On the other hand, by using a random quantizer argument it was shown~\cite{VBCI}
that by encoding vectors of infinite block length, it is possible to achieve distortions
\begin{eqnarray}
\lim_{R\rightarrow \infty}\bar{d}_0(R)2^{2R(1+a)} & = &
\frac{1}{4}\left(\frac{2^{2h(p)}}{2\pi e}\right), \nonumber \\
\lim_{R\rightarrow \infty}\bar{d}_s(R)2^{2R(1-a)} & = & \left(\frac{2^{2h(p)}}{2\pi e}
\right).     
\label{eqn-ecmdsq1}
\end{eqnarray}
Thus by using multiple description quantization it is possible to
simultaneously reduce the two-channel and side-channel granular distortions by 1.53 dB.

In single description quantization, an extra transmitted bit  reduces the squared error
distortion
by a factor of 4 (this is seen in (\ref{eqn11})). However, in multiple  description
quantization
there is additional flexibility. If each $R$ is increased by $1/2$ bit, the two-channel
distortion can be decreased
by $2^{-(1+a)}$ and the side distortion by $2^{-(1-a)}$, for any $a \in (0,1)$. This means
that 
by using an extra bit, the distortions $\bar{d}_0$ and $\bar{d}_s$ can be made to decrease by
different amounts
as long as the {\em product} decreases by a factor of 4.

\begin{figure}[htb]
\centerline{\psfig{figure=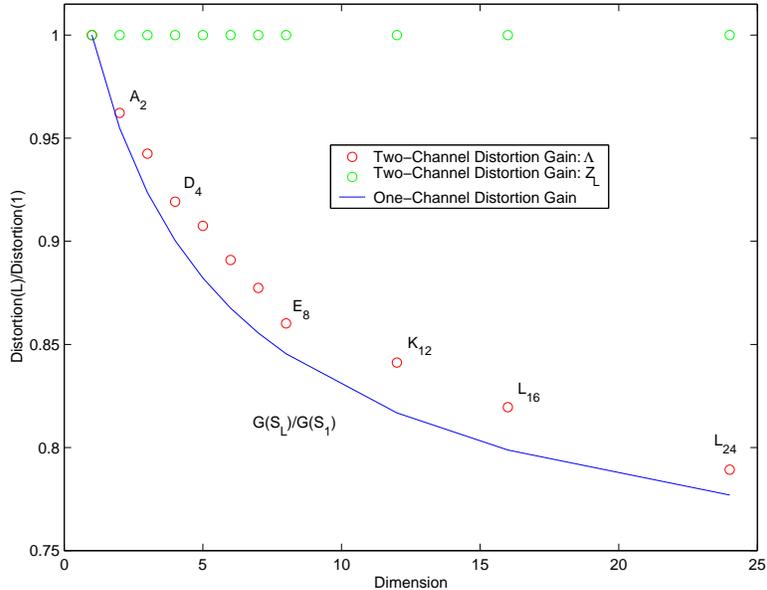,width=4in}}
\caption{\small Limiting two-channel and single-channel distortion ratios $\lim_{R\rightarrow \infty} \bar{d}_0(R,L)/\bar{d}_0(R,1)$ and
$\lim_{R\rightarrow \infty} \bar{d}_s(R,L)/\bar{d}_1(R,1)$ as a function of the lattice dimension $L$. 
The reduction in two-channel distortion is lattice dependent and is determined by the normalized
second moment of a Voronoi region of the lattice. The reduction in single-channel distortion
depends only on the dimension of the lattice.}
\label{figure-PackingGain}
\end{figure}

The goal of this paper is to give constructions for closing this ``1.53 dB'' gap and to
analyze the
resulting performance gains.
Our approach is as follows. From classical quantization theory, we know that the gap between scalar
quantization and the
rate distortion bound may be closed by using vector quantizers with lattice codebooks.
Certainly, by 
following this approach we can also close the gap between the two-channel distortion and the
rate-distortion
bound. In particular, this will allow us to replace the factor $(1/12)$ in the expression for
$\bar{d}_0$ in~(\ref{eqn-ecmdsq0})
with $G(\Lambda)$, the normalized second moment of the Voronoi region of a lattice point. The
main  question
we address here is that of simultaneously reducing $\bar{d}_s$. How can such a reduction be
achieved and 
what is the quantity that will replace the factor $(1/12)$ in the
expression for $\bar{d}_s$ in~(\ref{eqn-ecmdsq0})? We will show through a constructive
procedure that the
distortion $\bar{d}_s$ can be reduced by solving a specific labeling problem. To our
surprise, the quantity that 
replaces $(1/12)$ is $G(S_L)$,
the normalized second moment of a {\em sphere} in $L$ dimensions.

Let  $\bar{d}_0(R,L)$ and $\bar{d}_s(R,L)$ denote the two-channel and single channel
distortions at rate $R$ for an $L$-dimensional quantizer.  Fig.~\ref{figure-PackingGain} summarizes
the main results of the paper. In this figure we have plotted  the limit
of the normalized
two-channel and single-channel distortions, $\lim_{R\rightarrow \infty}\bar{d}_0(R,L)/\bar{d}_0(R,1)$ and 
$\lim_{R\rightarrow \infty}\bar{d}_s(R,L)/\bar{d}_s(R,1)$ respectively,
for lattices of various dimensions. It is seen that the limit for the
two-channel distortion is given by the ratio $G(\Lambda)/G(\ZZ)$,
which depends on the lattice, whereas for the side distortion the limit is $G(S_L)/G(S_1)$,
which is independent of the lattice.

\section{\bf Previous Work}
An achievable rate region for the multiple description problem was first given 
in~\cite{ECII} and it was shown in \cite{OZAI} that this region
coincides with the 
rate distortion region for a memoryless Gaussian source with a
squared-error distortion measure.
The problem of multiple description quantizer design, including a formulation and solution of
the underlying labeling
problem in one dimension, was presented in~\cite{Vai3}. An asymptotic performance analysis of
this quantizer
was presented in~\cite{BVI}.  A preliminary version of the work
presented here was first presented in~\cite{SerVaiSlo}. 

Lattice quantizers (for the single description problem) have been extensively studied.
In~\cite{Zador1} a random
quantization argument is used to give upper and lower bounds on the performance of quantizers for a
fixed
dimension. Detailed descriptions of the Voronoi regions of
specific lattices are given 
in~\cite{ConSlo1} and their second moments are evaluated. Fast quantizing algorithms for
lattice quantizers are given in~\cite{ConSlo2},~\cite{AdoulBarth:88},~\cite{Forney:88},~\cite{VardyBeery:93}. In~\cite{Gersho1}
it is conjectured that for
$r$th
power difference distortion measures, any optimal quantizer, in any dimension, has Voronoi
regions that are congruent to some polytope.  

Approaches to multiple description coding based on trellis coded quantization are
presented in~\cite{VBCI,JafTar:2}, those based on vector quantization
are presented in~\cite{FleEff:1} and approaches using forward
error correction are presented in~\cite{MohRisLad:1}.

There are also approaches to multiple description coding based on
subspace methods. One approach is to design a predictor or a
transform so as to achieve a correlation structure that allows
one half
of the prediction error samples or transform coefficients to be
predicted from the other half. Examples of this approach are
presented in \cite{IVI,OWVR,VSI}.

Another approach to multiple description  coding is based on
using overcomplete expansions. Here the idea is to construct a
redundant signal representation, in a way that allows the signal
to be estimated with a controlled amount of error when the
representation is incompletely
received~\cite{BalDauVai:1,ChouMehWang:1, GoyKovVet:1}.

\section{Organization of the Paper}
The multiple description vector quantizer is described, notation is established, certain regularity
assumptions for the labeling function are stated and preliminary expressions are
derived for rates and distortions in Section~\ref{section-preliminaries}. A detailed example with the
two-dimensional hexagonal lattice $A_2$ is presented along with some general theory in Section~\ref{section-A2-example}.
The necessary theory for general lattices is presented in Section~\ref{sec-gen}. An asymptotic analysis for
a fixed vector dimension is presented in Section~\ref{sec-asymp}. The paper is summarized along with conclusions and
directions for future work in Section~\ref{sec-summary}.

\section{Preliminaries}
\label{section-preliminaries}
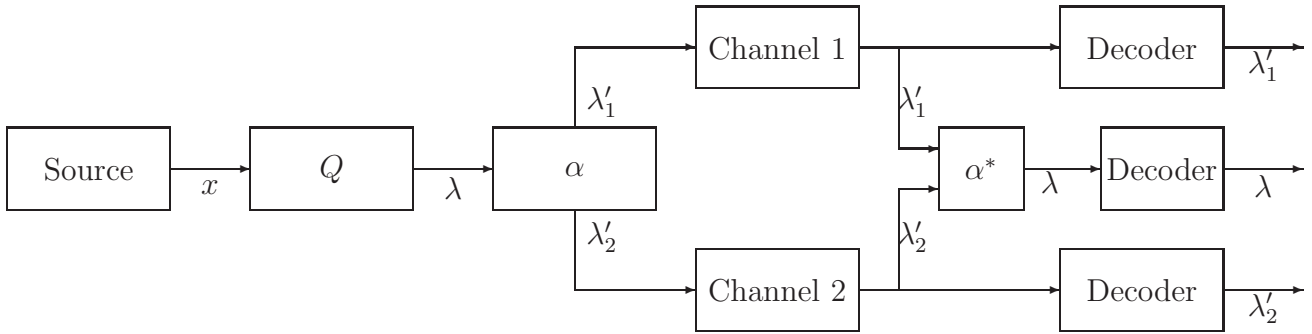
\begin{figure}[htb]
\setlength{\unitlength}{3347sp}%
\begingroup\makeatletter\ifx\SetFigFont\undefined
\def\x#1#2#3#4#5#6#7\relax{\def\x{#1#2#3#4#5#6}}%
\expandafter\x\fmtname xxxxxx\relax \def\y{splain}%
\ifx\x\y   
\gdef\SetFigFont#1#2#3{%
  \ifnum #1<17\tiny\else \ifnum #1<20\small\else
  \ifnum #1<24\normalsize\else \ifnum #1<29\large\else
  \ifnum #1<34\Large\else \ifnum #1<41\LARGE\else
     \huge\fi\fi\fi\fi\fi\fi
  \csname #3\endcsname}%
\else
\gdef\SetFigFont#1#2#3{\begingroup
  \count@#1\relax \ifnum 25<\count@\count@25\fi
  \def\x{\endgroup\@setsize\SetFigFont{#2pt}}%
  \expandafter\x
    \csname \romannumeral\the\count@ pt\expandafter\endcsname
    \csname @\romannumeral\the\count@ pt\endcsname
  \csname #3\endcsname}%
\fi
\fi\endgroup
\begin{picture}(9624,2424)(889,-2473)
\thinlines
\put(901,-1561){\framebox(1200,600){Source}}
\put(2701,-1561){\framebox(1200,600){$Q$}}
\put(4501,-1561){\framebox(1200,600){$\alpha$}}
\put(5301,-761){\makebox(0,0){$\lambda'_1$}}
\put(5301,-1761){\makebox(0,0){$\lambda'_2$}}
\put(6001,-661){\framebox(1200,600){\shortstack{Channel 1}}}
\put(6001,-2461){\framebox(1200,600){\shortstack{Channel 2}}}
\put(2101,-1261){\vector( 1, 0){600}}
\put(2401,-1400){\makebox(0,0){$x$}}
\put(3901,-1261){\vector( 1, 0){600}}
\put(4201,-1400){\makebox(0,0){$\lambda$}}
\put(5101,-961){\line( 0, 1){600}}
\put(5101,-361){\vector( 1, 0){900}}
\put(5101,-1561){\line( 0,-1){600}}
\put(5101,-2161){\vector( 1, 0){900}}
\put(8701,-661){\framebox(1200,600){Decoder}}
\put(8701,-2461){\framebox(1200,600){Decoder}}
\put(7801,-1561){\framebox(620,600){$\alpha^*$}}
\put(8621,-1361){\makebox(0,0){$\lambda$}}
\put(9001,-1561){\framebox(900,600){Decoder}}
\put(7201,-361){\vector( 1, 0){1500}}
\put(7601,-761){\makebox(0,0){$\lambda'_1$}}
\put(7601,-1761){\makebox(0,0){$\lambda'_2$}}
\put(7501,-361){\line( 0,-1){750}}
\put(7501,-1111){\vector( 1, 0){300}}
\put(7201,-2161){\vector( 1, 0){1500}}
\put(7501,-2161){\line( 0, 1){750}}
\put(7501,-1411){\vector( 1, 0){300}}
\put(8421,-1261){\vector( 1, 0){580}}
\put(9901,-361){\vector( 1, 0){600}}
\put(10201,-491){\makebox(0,0){$\lambda'_1$}}
\put(9901,-1261){\vector( 1, 0){600}}
\put(10201,-1391){\makebox(0,0){$\lambda$}}
\put(9901,-2161){\vector( 1, 0){600}}
\put(10201,-2291){\makebox(0,0){$\lambda'_2$}}
\end{picture}
\caption{Block diagram for a multiple description vector
quantizer.}
\label{fig-mdvq}
\end{figure}

A block diagram of a  multiple description vector quantizer (MDVQ) with a lattice codebook is shown in
Fig.~\ref{fig-mdvq}.

A source of information generates a sequence of independent, identically distributed
random variables with probability density function (pdf) $p$. This source is blocked off
into $L$-dimensional vectors $x=(x_1,x_2,\ldots,x_L)$. The $L$-fold pdf is denoted by $p_L$,
where $p_L(x)=\prod_{i=1}^Lp(x_i)$.  The vector  $x$  is quantized to the nearest vector
$\lambda$ in a lattice $\Lambda\subset \Reals^L$. We denote the quantizer mapping by $\lambda=Q(x)$.  Information about the selected code vector $\lambda$ is then
sent across the two channels, subject to rate constraints imposed by the individual channels.
This is done through a labeling function $\alpha$ followed by entropy coding.
The labeling function $\alpha$ maps $\lambda \in \Lambda$ to a pair
$(\lambda'_1,\lambda'_2) \in \Lambda' \times \Lambda'$, where $\Lambda'$ is a
sublattice of $\Lambda$ with index $N$. 
The component functions of $\alpha$ are denoted by $\alpha_1$ and $\alpha_2$, where
$\alpha_1(\lambda)=\lambda'_1$ and $\alpha_2(\lambda)=\lambda'_2$.
For simplicity we assume that
$\Lambda'$ is geometrically similar to $\Lambda$, i.e., $\Lambda'$ can be obtained by
scaling, rotating and possibly reflecting $\Lambda$.
Note that points in the lattice $\Lambda$ are denoted by $\lambda$,
possibly with subscripts, whereas
sublattice points will be denoted by $\lambda'$ or $\lambda''$, possibly with subscripts.

\begin{figure}[htb]
\centerline{\psfig{figure=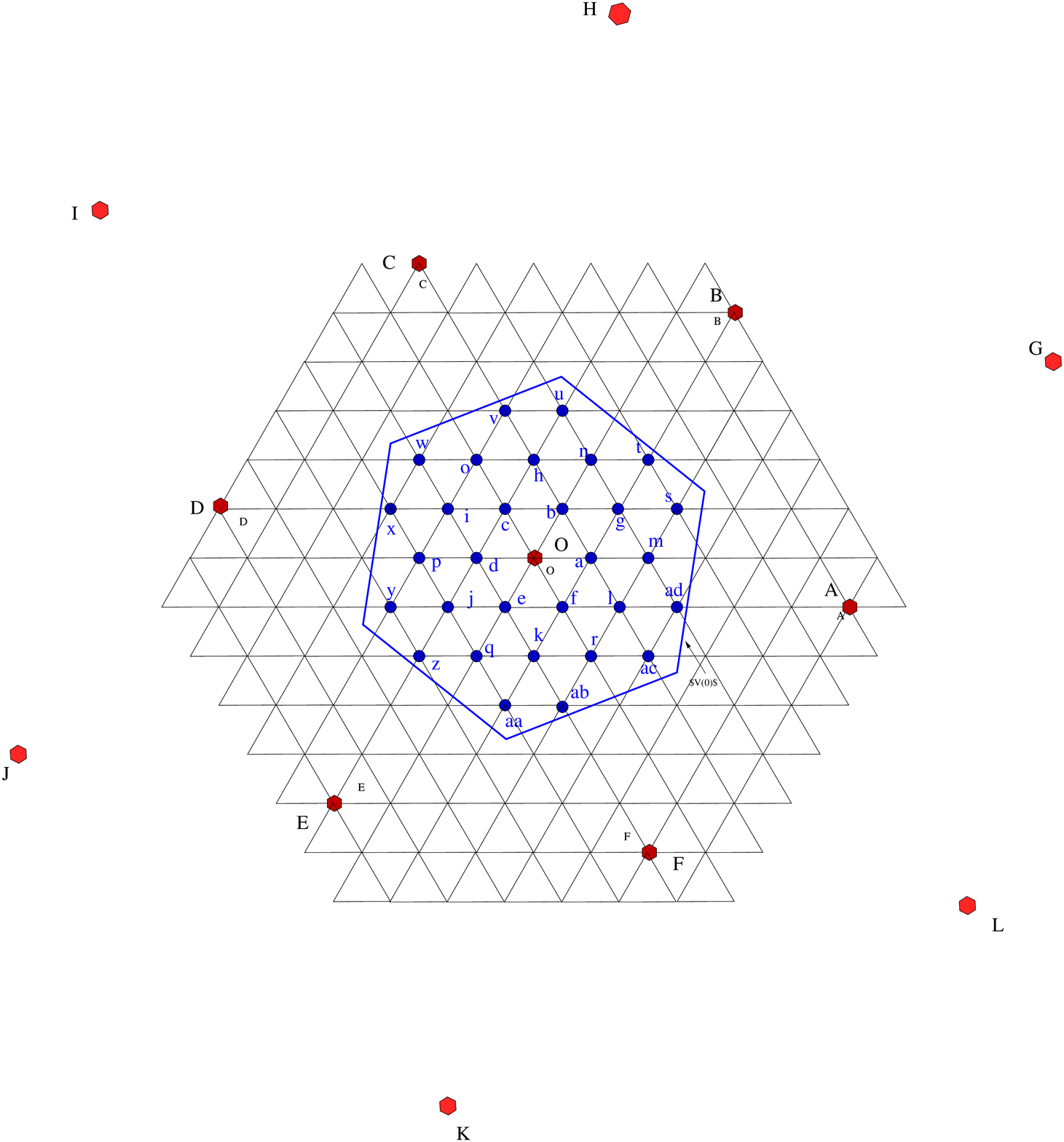,width=4in}}
\caption{This figure show a portion of the $A_2$ lattice (the points on
the grid line intersections), part of a geometrically similar sublattice of index
$31$, a discrete Voronoi set for a sublattice point (the lattice points in the hexagon) and the Voronoi set of a
sublattice point (all the points in the hexagon).}
\label{fig-notation1}
\end{figure}

In Fig.~\ref{fig-notation1}, a portion of the hexagonal lattice $A_2$ is illustrated, along with a geometrically
similar sublattice of index $31$. The lattice points lie at the intersection of the straight  lines
in the hexagonal grid (only some of the points are shown). The sublattice points are marked with upper-case
letters. Observe that the lattice is 31 times as dense as the sublattice, i.e., there are $31$ lattice
points for every sublattice point.

At the decoder, if only channel 1 works, the received information is used to decode 
$\lambda'_1$, and if only channel 2 works, the information received over channel
2 is used to decode $\lambda'_2$.  The mapping $\alpha$ is assumed
to be one-to-one so that if both channels work $\lambda$ can be recovered from
$(\lambda'_1,\lambda'_2)$.
(In practice, if only one channel is working it may be better to decode the
received vector to some function of $\la'_1$ or $\la'_2$ rather than to
$\la'_1$ or $\la'_2$ itself.
If $\la'_1$ is received but $\la'_2$ is not, for instance,
we would decode $\la'_1$ as the center of mass of all points
$\la \in \Lambda$ such that the first component of $\af ( \la )$
is $\la'_1$.
We will ignore this complication in order to simplify the analysis.)

Given $\Lambda$, $\Lambda'$ and $\alpha$, there are three distortions and two rates associated with an
MDVQ.
For a given $x$ mapped to the triple $(\lambda, \lambda'_1,\lambda'_2)$ by the MDVQ, the {\em two-channel
distortion} $d_0$ is given by $\|x-\lambda\|^2$, the channel 1 distortion $d_1$ by $\|x-\lambda'_1\|^2$
and the channel 2
distortion $d_2$ by $\|x-\lambda'_2\|^2$ (we assume that the inner product of $L$-dimensional vectors
$x=(x_1,x_2,\ldots,x_L)$ and $y=(y_1,y_2,\ldots,y_L)$ is given by $\langle x,y \rangle =(1/L)\sum_{i=1}^Lx_iy_i$ and the
corresponding norm is $\|x\|=\langle x,x \rangle ^{1/2}$, i.e., the inner product and norm are
dimension-normalized). The corresponding average distortions are denoted by $\bar{d}_0$,
$\bar{d}_1$ and $\bar{d}_2$.
We assume that an entropy coder is used in order to transmit the labeled
vectors at a rate arbitrarily close to the entropy, i.e., $R_i=\calH( \alpha_i(Q(\bX)) )/L$,
$i=1,2$, where $\calH(U)$, the entropy of the random variable $U$ taking values in alphabet ${\cal U}$ with
probability distribution $P$ is given by
$\calH(U)=-\sum_{u \in {\cal U}}P(u)\log P(u)$.

The problem is to design the labeling function $\alpha$
so as to minimize $\bar{d}_0$ subject to $\bar{d}_1\leq D_s$, $\bar{d}_2\leq D_s$,
and $R_i \leq R$, $i=1,2$, for specified values of the rate $R$ and  distortion $D_s$.

The Voronoi (or nearest neighbor) region $V(\lambda)$ of a point $\lambda \in \Lambda$ is defined
to be
\begin{equation}
V(\lambda)\stackrel{def}{=}\{x~:~\|x-\lambda\|\leq \|x-\tilde\lambda\|,~\forall \tilde\lambda \in \Lambda\}.
\label{eqn-app-1}
\end{equation}
With each sublattice point $\lambda' \in \Lambda'$, we associate a discrete Voronoi set (with $N$ elements)
\begin{equation}
V_0(\lambda')\stackrel{def}{=}\{\lambda\in \Lambda~:~\|\lambda-\lambda'\| \leq \|\lambda-\lambda'' \|,\forall \lambda'' \in \Lambda' \}.
\label{eqn-app-1a}
\end{equation}
In (\ref{eqn-app-1}) and (\ref{eqn-app-1a}) ties (i.e., points for which equality holds in the defining condition) 
are broken in some prearranged manner. 
The existence of sublattices for which no ties occur is discussed in~\cite{ConRaiSlo1}.
The Voronoi region $V(0)$ and the discrete Voronoi region $V_0(0)$ are both illustrated in Fig.~\ref{fig-notation1}.

We regard the label for $\la$ as a directed edge $\oer = (\alpha_1 (\la ), \alpha_2 ( \la ))$ of the graph with vertex set $\La'$.
The corresponding unordered pair $e=\{\alpha_1(\lambda),\alpha_2(\lambda)\}$ will be referred to as the 
undirected edge or undirected label associated with $\lambda$. The essential difference between a directed edge or label $(\lambda_a',\lambda_b')$ and the undirected edge $\{\lambda_a',\lambda_b'\}$ is that for the directed edge there is an implicit association between 
edge component and channel ($\lambda_a'$ is sent
on channel 1 and $\lambda'_b$ is sent on channel 2) whereas for the undirected edge no association is implied. 
Graphically, an edge connecting two sublattice points $\lambda'_a$ and
$\lambda'_b$, with an arrow pointing from $\lambda'_a$ to $\lambda'_b$, indicates that $\lambda'_a$ is sent on channel 1 and $\lambda'_b$ is sent
on channel 2, or equivalently, the directed edge is $(\lambda'_a,\lambda'_b)$.
The two directed versions of an (undirected) edge $e$
will be denoted $\oer$ and $\oel$.

For a given labeling function $\alpha$,
an associated undirected edge labeling function $\alpha_u$ is defined as follows:
if $\alpha(\lambda)=\oer$, then $\alpha_u(\lambda)=e$, i.e., $\alpha_u$ maps $\lambda$ to its undirected label.
Note that if $\alpha(\lambda_1)=\oer$ and $\alpha(\lambda_2)=\oel$, then
$\alpha_u(\lambda_1)=\alpha_u(\lambda_2)=e$, i.e., $\alpha_u$ is not one-to-one.
A directed edge is uniquely associated with 
a lattice point whereas an undirected edge will in general be associated with two lattice points, one for
each orientation of the edge. The main reason for introducing $\alpha_u$ is that its construction logically precedes
that of $\alpha$.

\begin{figure}[htb]
\hspace*{\fill}\psfig{figure=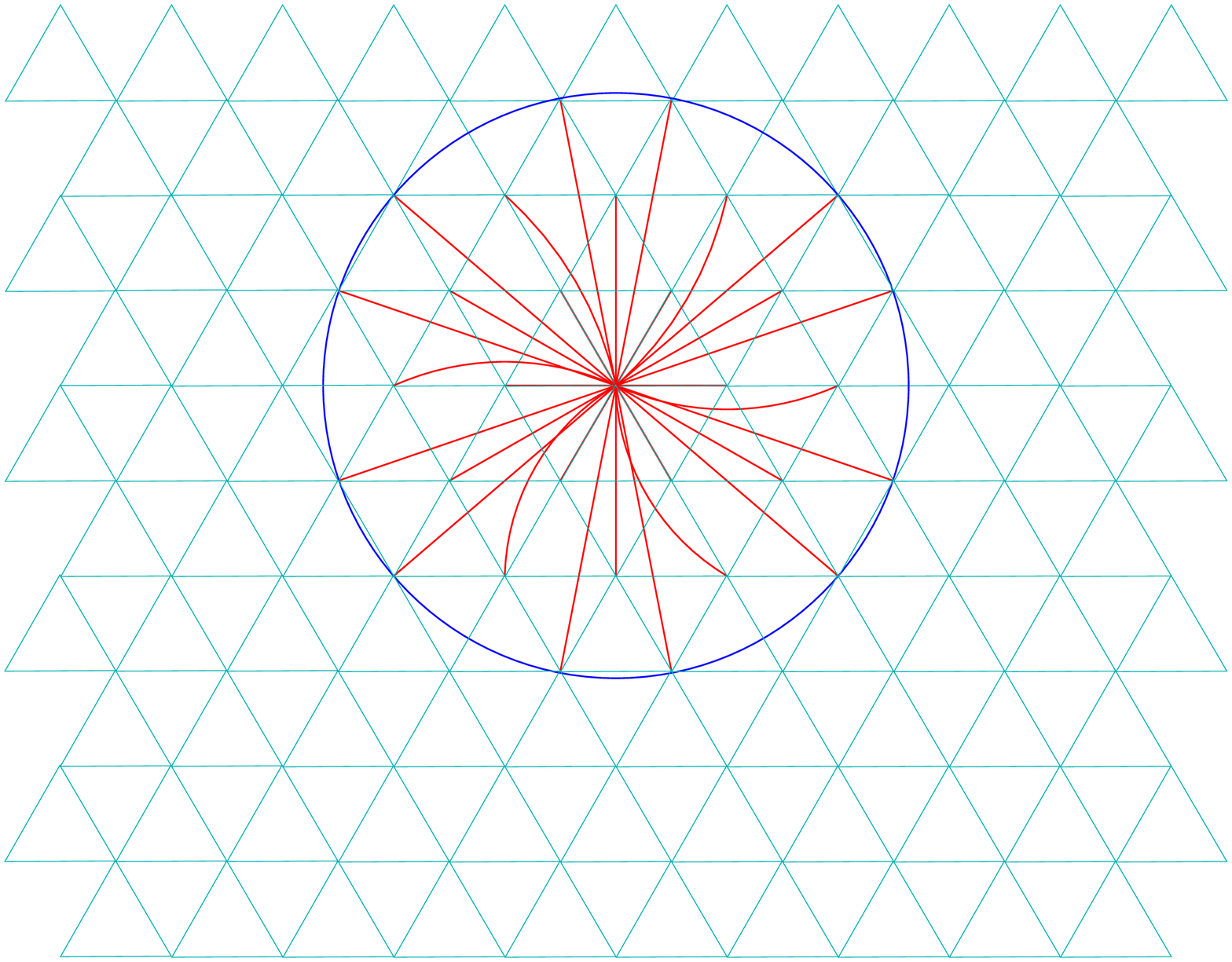,width=3in}\hspace{\fill}\psfig{figure=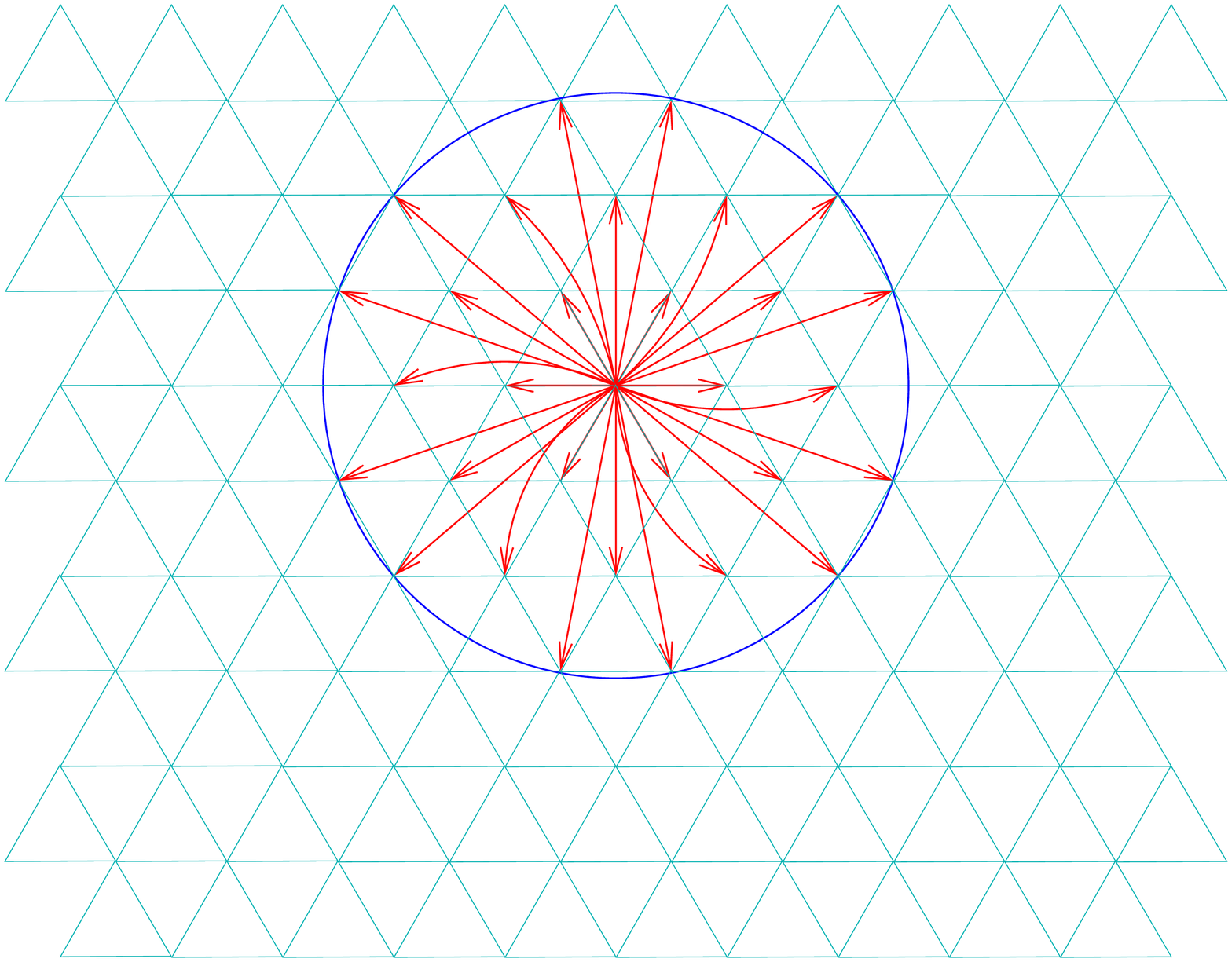,width=3in} \hspace*{\fill}

\caption{Various edge sets are illustrated for a sublattice of $A_2$ with $N=31$. For clarity only the sublattice points are shown. The set $\calE_{u|\alpha_1=0}$ is shown on the left, 
$\calE_{d|\alpha_1=0}$ on the right. The origin is located at the center of the circle.}
\label{fig-edgeset-1-2}
\end{figure}

Certain sets associated with the two maps will play a central role in the development. The first is
$\calE_d(\lambda')=\alpha(V_0(\lambda'))$, the set of all labels for points in $V_0(\lambda')$
(the subscript $d$ indicates that the edges are directed). The corresponding set of undirected edges is denoted by
$\calE_u(\lambda')$. More specifically, $\calE_u(\lambda')=\alpha_u(V_0(\lambda'))$.
The set of all (directed) labels is denoted by $\calE_d
=\bigcup_{\lambda' \in \Lambda'}\calE_d(\lambda')$ and the set of
all undirected labels by $\calE_u=\bigcup_{\lambda' \in
\Lambda'}\calE_u(\lambda')$. It is also useful to define the
restriction $\calE_d|_{\alpha_1=0}=\{(\lambda'_1,\lambda'_2) \in
\calE_d,~\lambda'_1=0\}$. The set of undirected edges in
$\calE_d|_{\alpha_1=0}$ will be denoted by
$\calE_u|_{\alpha_1=0}$. The sets ${\cal E}_u|_{\alpha_1=0}$ and ${\cal E}_{d}|_{\alpha_1=0}$ are illustrated in
Fig.~\ref{fig-edgeset-1-2}.

The {\em reuse index} associated with a label $\lambda'$
and channel $i$ is defined to be $N_i(\lambda')=|\{\lambda~:~\alpha_i(\lambda)=\lambda'\}|$,
the number of lattice points $\lambda$ for which
$\alpha_i(\lambda)=\lambda'$.

In order to render the problem tractable, we assume that the 
labeling function has the following properties. 
\begin{itemize}
\item{\bf Property 1}:
The reuse index $N_i(\lambda')=N$, for all $\lambda' \in
\Lambda'$ and $i=1,2$. In other words each channel $i$ label is reused exactly $N$ times.
\item{\bf Property 2}: The shift property: $\alpha_u(\lambda+\lambda')=\alpha_u(\lambda)+\lambda'$, for all
$\lambda \in \Lambda$ and $\lambda' \in \Lambda'$.
\item{\bf Property 3}: Each undirected edge $e=\{\lambda'_1,\lambda'_2\}$, $\lambda_1'\neq \lambda_2'$, labels two points, $\lambda_a$
and
$\lambda_b$, and $\lambda_a+\lambda_b=\lambda'_1+\lambda'_2$. 
\end{itemize}
The first assumption makes it easy to parametrize the tradeoff between
the side and central distortions.
The second assumption reduces the labeling problem
to one of labeling a finite set.
The third assumption is
a simple way to achieve exact balance between the two descriptions (its
implications will become clearer in Section~\ref{sec-balance}).

It will prove useful to have the following definitions of equivalence.

\noindent
{\bf Definition 1}: Two lattice points $\lambda_1$ and
$\lambda_2$ are said to be equivalent if either $\lambda_1$ and $\lambda_2$
or $\lambda_1$ and $-\lambda_2$ lie in the same
coset of $\Lambda$ relative to the sublattice $\Lambda'$.

\noindent
{\bf Definition 2}: Two sublattice edges $e_1$ and $e_2$ are
said to be equivalent if they are parallel and of equal length,
or equivalently if $e_1 +\lambda'=e_2$ for some $\lambda' \in
\Lambda'$.

\noindent
The equivalence class of an object will be indicated by square
brackets. Thus
$[\lambda]$ is the equivalence class of $\lambda$ and $[e]$ is
the equivalence class of $e$, or equivalently, $[\lambda]
\in \Lambda / \Lambda'$ and $[e]\in\calE_u / \Lambda'$.

\subsection{\bf Distortion Computation}

The average two-channel distortion $\bar{d}_0$ is given by
\begin{equation}
\bar{d}_0=\sum_{\lambda \in
\Lambda}\int_{V(\lambda)}\|x-\lambda\|^2 p_L(x) dx.
\label{eqn-app-2}
\end{equation}
Since the codebook of the quantizer is a lattice, all the Voronoi regions are congruent. Furthermore,
upon assuming that each Voronoi region is small so that $p_L(x)\approx p_L(\lambda)$ for $x \in V(\lambda)$ and upon letting $\nu$ denote the $L$-dimensional
volume of a Voronoi region, we obtain the following expression for
the two-channel distortion \cite{Gersho1}
\begin{equation}
\bar{d}_0\approx\frac{\int_{V(0)}\|x\|^2dx}{\nu} ~,
\end{equation}
which in terms of the normalized second moment $G(\Lambda)$, defined by
\begin{equation}
G(\Lambda)\stackrel{def}{=}\frac{\int_{V(0)}\|x\|^2 dx}{\nu^{1+2/L}},
\label{eqn-cen-dis}
\end{equation}
is given by
\begin{equation}
\bar{d}_0\approx G(\Lambda)\nu^{2/L}.
\label{eqn-d0}
\end{equation}

We now derive expressions for the average distortions $\bar{d}_1$ and $\bar{d}_2$ and
$\bar{d}_s\stackrel{def}{=}(\bar{d}_1+\bar{d}_2)/2$.
When only description 1 is available, the distortion is given by
\begin{eqnarray}
\bar{d}_1 & = & \sum_{\lambda \in \Lambda}
\int_{V(\lambda)}\|x-\alpha_1(\lambda)\|^2 p_L(x) dx \nonumber \\
& = & \sum_{\lambda \in \Lambda}
\int_{V(\lambda)}\|x-\lambda+\lambda-\alpha_1(\lambda)\|^2 p_L(x) dx \nonumber \\
& = & \sum_{\lambda \in
\Lambda}\int_{V(\lambda)}\|x-\lambda\|^2p_L(x) dx +   \sum_{\lambda \in \Lambda}
\int_{V(\lambda)}\|\lambda-\alpha_1(\lambda)\|^2 p_L(x) dx +\\
& &   2\sum_{\lambda \in \Lambda} \int_{V(\lambda)}\left\langle
x-\lambda,\lambda-\alpha_1(\lambda) \right\rangle p_L(x) dx \nonumber \\
& = & \bar{d}_0 + \sum_{\lambda \in \Lambda} \|\lambda-\alpha_1(\lambda)\|^2 P(\lambda) + 2\sum_{\lambda
\in \Lambda}
\left\langle \int_{V(\lambda)}xp_L(x)dx-\int_{V(\lambda)}\lambda
p_L(x) dx, \lambda-\alpha_1(\lambda) \right\rangle \nonumber \\
& \stackrel{(a)}{=} & \bar{d}_0 + \sum_{\lambda \in \Lambda} \|\lambda-\alpha_1(\lambda) \|^2 P(\lambda),
\label{eqn-app-3}
\end{eqnarray}
where $P(\lambda)=Pr(Q(X)=\lambda)$, and (a) follows by assuming that
$\lambda$ is the
{\em centroid} of its Voronoi region. This is true for the uniform density. For nonuniform densities,
there is an 
error term which goes to zero with the size of the Voronoi region.
The first term in (\ref{eqn-app-3}) is the two-channel distortion and the second term is the excess
distortion
which is incurred when channel 2 fails. Note that for a given $\Lambda$, only the excess distortion term
is affected
by the labeling function $\alpha$. 
From (\ref{eqn-app-3}), it follows that 
\begin{eqnarray}
\bar{d}_s & = & \bar{d}_0 + (1/2) \sum_{\lambda \in \Lambda} (\|\lambda -\alpha_1(\lambda)\|^2 +\|\lambda -\alpha_2(\lambda)\|^2 )P(\lambda)\\
 & = & \bar{d}_0 + (1/2) \sum_{\lambda' \in \Lambda'} \sum_{\lambda \in V_0(\lambda')}(\|\lambda -\alpha_1(\lambda)\|^2 +\|\lambda -\alpha_2(\lambda)\|^2 )P(\lambda).
\label{eqn-app-4}
\end{eqnarray}

We also introduce the following notation: $d_1(\lambda,\oer)=\|\lambda-\lambda_1'\|^2$ and $d_2(\lambda,\oer)=\|\lambda-\lambda_2'\|^2$ and $d_s(\lambda,\oer)=(d_1(\oer)+d_2(\oer))/2$, where $\oer=(\lambda_1',\lambda_2')$. Note that $d_s(\lambda,\oer)=d_s(\lambda,\oel)$. Hence when appropriate, we will write $d_s(\lambda,e)$. Also, when the edge associated with the lattice point is
clear from the context, we will write $d_s(\lambda)$ or $d_s(\oer)$ instead of $d_s(\lambda,\oer)$.
It is useful (as a design guide and for the asymptotics which follow) to write down a slightly different expression for the
side distortion where the sum is taken over the edge set $\calE_d=\bigcup_{\lambda'\in \Lambda'}\calE_d(\lambda')$:
\begin{eqnarray}
\bar{d}_s & {=} & \bar{d}_0+\sum_{\oer \in \calE_d}d_s(\oer)P(\oer)\nonumber \\
 & {=} & \bar{d}_0+\sum_{\lambda' \in \Lambda'}\sum_{\oer \in \calE_d(\lambda')}d_s(\oer)P(\oer)
\label{eqn-app-5}
\end{eqnarray}
where $P(\oer)$ is equal to the probability of the lattice point that the edge labels, i.e.,
 $P(\oer)=Pr(Q(X)=\lambda)$.

\subsection{\bf Rate Computation}

Expressions for the rate (in bits/sample) will be derived next. Let $R_0$ be
the rate required to address the
two-channel codebook for a single channel system\footnote{This quantity is
useful for evaluating the 
two-channel distortion as well as for evaluating the rate overhead associated
with the multiple description
scheme.}.
We will first derive an expression for $R_0$ and then determine the
(per-channel) rate $R$ of the multiple
description system. 

In order to derive expressions for $R_0$ and $R$, we use the fact that each
quantizer bin
has identical volume $\nu$ and that $p_L(x)$ is approximately constant over
Voronoi regions of
the sublattice $V(\lambda')$.
The second assumption is valid in the limit as the Voronoi regions become
small and is standard
in asymptotic quantization theory.

The rate $R_0=\calH(Q(X))$ is given by \cite{Gray:90}
\begin{eqnarray}
R_0 & = & -(1/L)\sum_{\lambda}\int_{V(\lambda)}p_L(x)dx \log_2
\int_{V(\lambda)} p_L(x) dx \nonumber \\
  & \approx & -(1/L) \sum_{\lambda}\int_{V(\lambda)}p_L(x)dx \log_2
p_L(\lambda) \nu \nonumber \\
  & \approx & h(p)-(1/L)log_2(\nu).
\label{eqn-rate-0}
\end{eqnarray}

For $R$, we evaluate the entropy $\calH(\alpha_1(Q(X)))$ and then use the
approximation that
$p_L(x)$ is roughly constant over each Voronoi region of $\Lambda'$ to get
\begin{eqnarray}
R & = & -(1/L)\sum_{\lambda'\in \Lambda'}\left(\sum_{\lambda \in
\alpha_1^{-1}(\lambda')}\int_{V(\lambda)}p_L(x)dx
\right)\log_2\left[\left(\sum_{\lambda \in
\alpha_1^{-1}(\lambda')}\int_{V(\lambda)}p_L(x)dx \right) \right] \nonumber
\\
  & \approx & -(1/L)\sum_{\lambda'\in \Lambda'}\left(\sum_{\lambda \in
\alpha_1^{-1}(\lambda')}\int_{V(\lambda)}p_L(x)dx
\right)\log_2\left(p_L(\lambda')N\nu \right)
\nonumber \\
& \approx & -(1/L)\sum_{\lambda'\in \Lambda'}\left(\sum_{\lambda \in
\alpha_1^{-1}(\lambda')}\int_{V(\lambda)}p_L(x)\log_2(p_L(\lambda'))dx\right) -
(1/L)\log_2(N\nu)
\nonumber \\
& \approx & h(p)-(1/L)\log_2(N\nu).
\label{eqn-rate-1}
\end{eqnarray}
Observe that in the above equation, the term $N\nu$ is simply the volume of a
fundamental
region for the sublattice $\Lambda'$ (since it has index $N$ in $\Lambda$).
Upon writing (\ref{eqn-rate-1}) in terms of $R_0$ we obtain
\begin{equation}
R=R_0-(1/L)\log_2(N).
\end{equation}
A single channel system would have used $R_0$ bits/sample to achieve the same $\bar{d}_0$. Instead a multiple
description system uses a
total of $2R=2R_0-(2/L)\log_2(N)$ bits/sample, and so
the rate overhead is $R_0-(2/L)\log_2(N)$.

\section{A Labeling Function for $A_2$}
\label{section-A2-example}

We now look for a labeling function $\alpha$ for which $\sum_{\overrightarrow{e} \in \calE_d(\lambda')}d_s(\oer)$ 
is minimized and is independent of $\lambda'$.  Since the details are complicated, we will 
work out the first example--for the 
hexagonal lattice $A_2$ --quite carefully. We will then identify certain
general principles and use them to construct labelings for other lattices. 

The lattice $A_2$ may be considered to be a subset of $\Reals^2$ or as a subset of $\Complexes$. Since each
approach has its advantages, we will switch back and forth between the two representations.
We consider the lattice $A_2$ at unit scale to be generated by the vectors $\{1,\omega\}\subset \Complexes$, where
$\omega=-1/2+i\sqrt{3}/2$. The associated Gram matrix is $\left(\begin{smallmatrix}1 & -1/2 \\ -1/2 & 1 \end{smallmatrix}\right)$ and
the fundamental volume is $\sqrt{3}/2$.
A sublattice $\Lambda'$ of a lattice $\Lambda$ is said to be {\em geometrically similar} to $\Lambda$
if it can be obtained by scaling and rotating and/or reflecting $\Lambda$ \cite{SPLAG}.
To be more precise, if a matrix $G'$ generates
$\Lambda'$ and $G$ generates $\Lambda$, then $\Lambda'$ is geometrically similar to $\Lambda$ if and only if
$G'=cUGB$, for some nonzero scalar $c$, integer matrix $U$ with determinant $\pm 1$, and real orthogonal
matrix $B$. The index $N$ is defined as the ratio of the fundamental volumes of $\Lambda'$ and $\Lambda$ and is given in 
terms of the scale factor $c$ by  $N=c^2$.
It can be shown~\cite{BerWriSlo} that $\Lambda'$ is similar to $\Lambda$ if and only if $N$ is of the form
$a^2-ab+b^2$, $a,b \in \ZZ$; if this holds then $\Lambda'$ is generated by $\bu=a+b\omega$ and $\bv=\omega(a+b\omega)$.
In addition to this restriction on $N$, we will require, for convenience only, that $N=\sum_{i=0}^K A_i$, where
$A_i$ is the number of lattice points at squared distance $i$ from the origin. In other words, we require that
$N$ is the number of points in the first $K$ shells of the lattice,
for some $K=K(N)$.
There are heuristic arguments, to be presented elsewhere, which suggest that there are infinitely many values of $N$ with this property.
For example, $N=31$ has this property, since
$N= A_0 + A_1 +A_2+A_3 +A_4 = 1+ 6+6+6+12 = 31$, and 31 is
also of the form $a^2-ab+b^2$, with
$a=5, b=-1$.

\subsection{\bf An Example}
\label{sec-example-a2}
\begin{figure}[htb]
\centerline{\psfig{figure=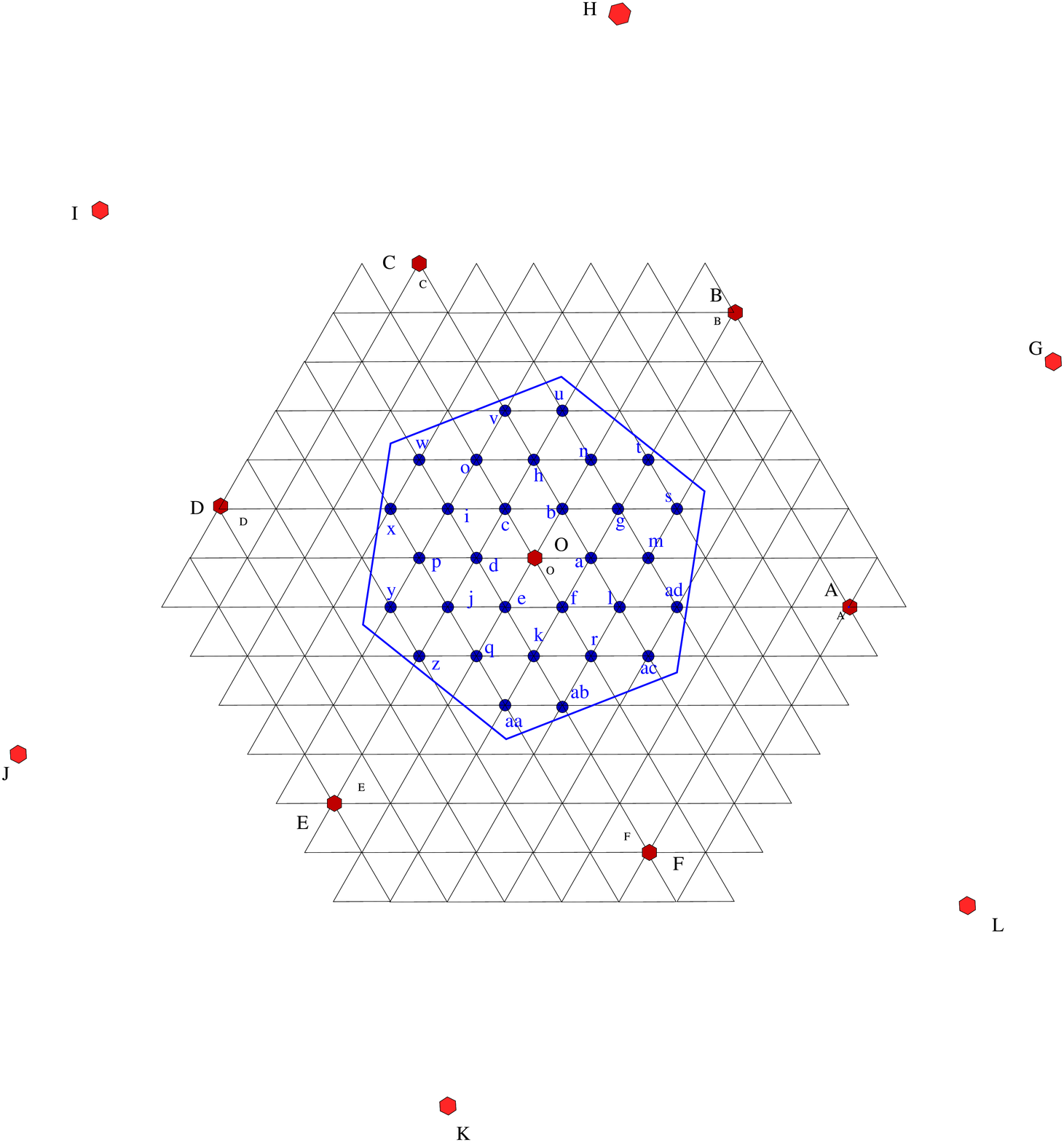,width=4.5in}}
\caption{Example to illustrate the mechanics of the labeling for the hexagonal lattice $A_2$.}
\label{fig-example-1}
\end{figure}

We now present an example of a labeling function with a reuse index $N=31$. 
The following are the steps in constructing the labeling function.
\begin{enumerate}
\item Find a sublattice with index equal to the reuse index.
\item Determine the discrete Voronoi set $V_0(0)$.
\item Determine an undirected label for every point $\lambda \in V_0(0)$.
\item Extend the labeling to the entire lattice using the shift property of the undirected
labels.
\item Given $\lambda$ and its undirected label $e$, determine the correct directed label $\oer$ (i.e., determine which endpoint of $e$ is to be sent on channel 1 and which endpoint on channel 2).
\end{enumerate}
A sublattice $\Lambda'$ of index equal
to $31$ may be obtained by considering all points of the form $a\bu + b \bv$, with $\bu=5-\omega$ and $\bv=\omega(5-\omega)$.
A portion of this lattice together with some sublattice points is shown in Fig.~\ref{fig-example-1}.
Lattice points have been labeled with lower-case letters
$a,b,c,\ldots$ and sublattice points with
upper-case letters $O, A, B, C, \ldots $.
In order to fix the coordinate system, note that the sublattice point $O$ is the origin of the complex plane and the lattice
points $a$ and $c$ have representations $1$ and
$\omega$ respectively. Relative to this basis, the representations of sublattice points $A$ and $C$ are $5-\omega$ and
$\omega(5-\omega)$, respectively. 
The discrete Voronoi set $V_0(0)= \{O,a,b,c,\ldots,y,z,aa,ab,ac,ad\}$ is also shown in Fig.~\ref{fig-example-1}.
Note that $|V_0(0)|=31$.
\begin{figure}[htb]
\begin{center}
\input{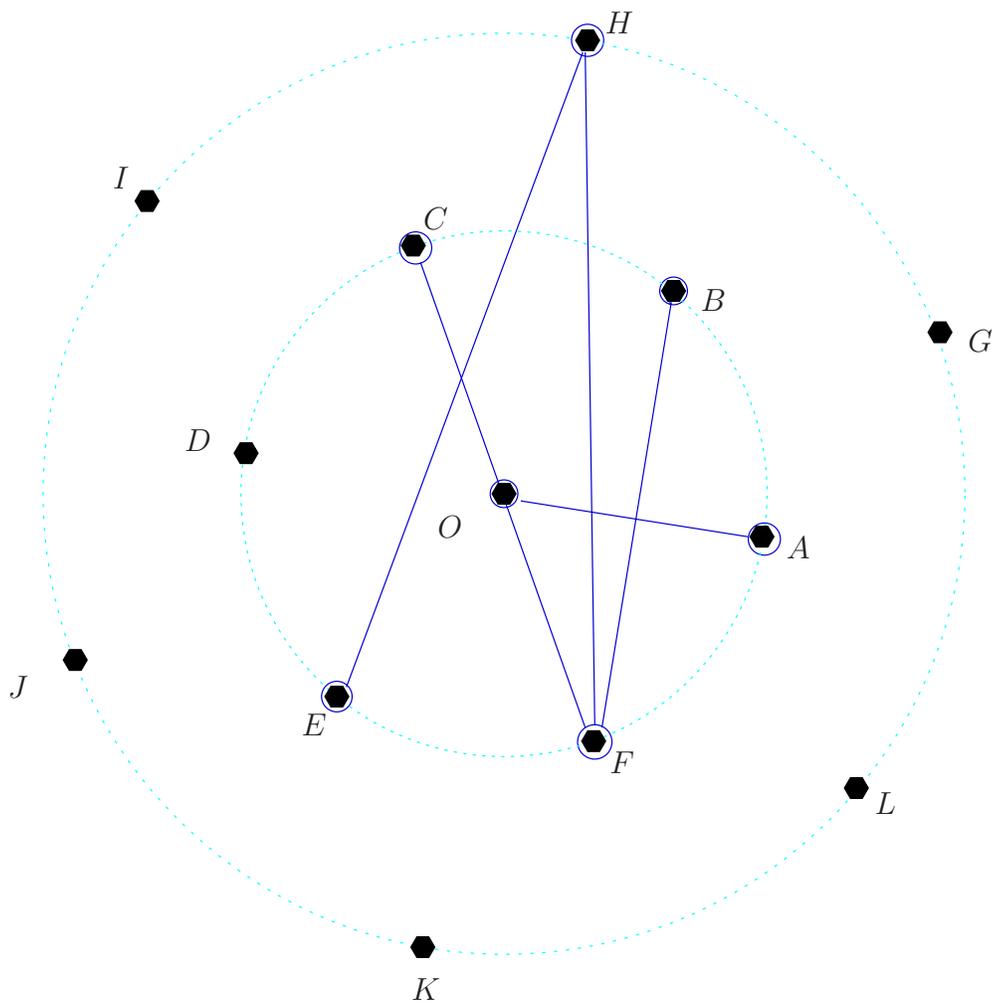}
\end{center}
\caption{The set of undirected edges ${\cal E}_u(0)$ for the example in Section~\ref{sec-example-a2}.
For clarity only a subset of ${\cal E}_u(0)$ is shown. The entire set may be obtained by
rotating each edge around the origin (sublattice point $O$) by multiples of $\pi/3$ radians. 
${\cal E}_u(0)$ consists of 28 distinct edges obtained by counting six edges for each edge in the figure,
except for the diameter, for which only three edges are counted.
The edge $\{O,O\}$ is not shown.}
\label{fig-EdgeSet}
\end{figure}
Points in $V_0(0)$ will be labeled using {\em directed} edges obtained from the following
set of $28$ {\em undirected} edges:
\begin{eqnarray}
\calE_u(0) & = & \{\{O,O\},\\
           &   & \{O,A\},\{O,B\},\{O,C\},\{O,D\},\{O,E\},\{O,F\}, \nonumber \\
           &   &   \{A,C\},\{B,D\},\{C,E\},\{D,F\},\{E,A\},\{F,B\}, \nonumber \\
           &   &   \{A,D\}, \{B,E\},\{C,F\}, \nonumber  \\
           &   &   \{G,D\},\{G,E\},\{H,E\},\{H,F\},\{I,F\}, \{I,A\}, \nonumber \\
           &   &   \{J,B\},\{J,A\},\{K,C\},\{K,B\},\{L,C\},\{L,D\}\},
\end{eqnarray}
which are illustrated in Fig.~\ref{fig-EdgeSet}. Three edges in
this set, namely $\{A,D\}$, $\{B,E\}$ and $\{C,F\}$, will be used twice (i.e., both orientations will be
used) in order to obtain $31$ labels.

\begin{table}[htb]
\centering
\begin{tabular}{|c|c|c|c||c|c|c|c|}\hline
Lattice & Label        & Color & Label &   Lattice & Label        & Color & Label  \\
Point   & (undirected) &       & (directed) &   Point   & (undirected) & & (directed) \\
 $\lambda$     & $e=\alpha_u(\lambda)$& $c(e)$& $\oer=\alpha(\lambda)$       & $\lambda$       & $e=\alpha_u(\lambda)$          & $c(e)$& $\oer=\alpha(\lambda)$        \\
                                                                                                    \hline \hline
 O             &  \{O,O\}     &       & (O,O) &                 &              &       &        \\
 a             &  \{O,A\}     & 0     & (O,A) & d               & \{O,D\}      & 1     & (D,O)  \\
 b             &  \{O,B\}     & 0     & (O,B) & e               & \{O,E\}      & 1     & (E,O)  \\
 c             &  \{O,C\}     & 0     & (O,C) & f               & \{O,F\}      & 1     & (F,O)  \\
 g             &  \{A,C\}     & 0     & (A,C) & j               & \{D,F\}      & 1     & (F,D)  \\
 h             &  \{B,D\}     & 0     & (B,D) & k               & \{E,A\}      & 1     & (A,E)  \\
 l             &  \{F,B\}     & 0     & (F,B) & i               & \{C,E\}      & 1     & (E,C)  \\
 n             &  \{A,D\}     & 0     & (A,D) & q               & \{A,D\}      & 0     & (D,A)  \\
 o             &  \{B,E\}     & 0     & (B,E) & r               & \{B,E\}      & 0     & (E,B)  \\
 p             &  \{C,F\}     & 0     & (C,F) & m               & \{C,F\}      & 0     & (F,C)  \\
 x             &  \{A,I\}     & 0     & (I,A) & ad              & \{D,L\}      & 1     & (D,L)  \\
 s             & \{D,G\}      & 0     & (G,D) & y               &  \{A,J\}     & 1     & (A,J)  \\
 t             & \{E,G\}      & 0     & (G,E) & z               &  \{B,J\}     & 1     & (B,J) \\
 u             & \{E,H\}      & 0     & (H,E) & aa              &  \{B,K\}     & 1     & (B,K) \\
 v             & \{F,H\}      & 0     & (H,F) & ab              &  \{C,K\}     & 1     & (C,K) \\
 ac            &  \{C,L\}     & 0     & (L,C) & w               & \{F,I\}      & 1     & (F,I)   \\ \hline \hline

\end{tabular}
\caption{
Labels for points in the discrete Voronoi set $V_0 (0)$.
The first step is to assign an undirected label to each lattice point in $V_0(0)$. In this case
the labeling was done by hand, though an optimal procedure is described in Section~\ref{sec-alphau}.
It is important that if two lattice points sum to zero, the midpoints of their labels also
sum to zero (equivalent points are shown in the same row). The second step is to assign a color $c(e)$
($0$ or $1$) to each edge and then to use this to determine the orientation of each edge. If $c(e)=0$,
the endpoint of $e$ which is closer to $\lambda$ becomes the channel 1 label and the
endpoint which is farther from $\lambda$ becomes the channel 2 label.
On the other hand if $c(e)=1$, the endpoint of $e$ which is closer becomes the channel 2 label and that
which is farther becomes the channel 1 label. The labeling is extended to the entire lattice using
the shift property of $\alpha_u$.
}
\label{table-alphau}
\end{table}
\clearpage
Each point in $V_0(0)$ is then associated with an undirected edge in $\calE_u(0)$ as shown in
Table~\ref{table-alphau}. This association was done by hand (but see Section 4.2.2).
The only constraint is that
equivalent points are mapped to equivalent edges--this can be clearly seen in the table, 
where equivalent lattice points have been placed in the same row.
In order to complete the labeling we need to assign a {\em directed}
edge to each lattice point based on the color $c(e)$ ($0$ or $1$ valued) associated with that edge.
This is illustrated for two cases--for a point that lies in $V_0(0)$ and for a point that lies
outside $V_0(0)$.

First consider the point $ac=1-2\omega$ which has been assigned the edge $\{C,L\}=\{1+6\omega, 4-7\omega\}$.
This edge has midpoint $\mu=5/2-1/2\omega$. The color of the edge $\{C,L\}$ is obtained in terms of the first component of the
midpoint of the edge ($5/2$) and the difference between the first components of the edge endpoints
($4-1=3$) by determining whether $\lfloor(5/2)/3\rfloor$ is odd or even. In this case it is even, hence
it is assigned the color $0$. (If it had been odd, it would have been assigned the color $1$.
An expression for the coloring rule is given later in (\ref{eqn-ColoringRule}).)
The specific orientation of this edge is then obtained by determining which endpoint of $\{C,L\}$ is
closer to $ac$. In this case $L=4-7\omega$ is closer to $ac$ than $C=1+6\omega$. Since the edge has color equal to $0$,
the closer endpoint becomes the channel 1 label, i.e., $ac$ is assigned the edge $(L,C)$ (this is the edge direction rule  (\ref{eqn-DirectionRuleA}), which is explained in detail in Section~\ref{section-edge-direction}). 

Now consider the lattice point $\lambda=18+10\omega \not \in V_0(0)$. The nearest lattice point is $\lambda'=17 +9 \omega$.
We compute the difference
$\lambda-\lambda'=(1+\omega)$ (the point $b$) and use this to  look up the corresponding undirected label in Table~\ref{table-alphau}.
This gives $\alpha_u(1+\omega)=\{O,B\}=\{0,6+5\omega\}$. Using the shift property, the undirected label for $\lambda$
is then obtained by  shifting the edge by $\lambda'$  to give $\alpha_u(\lambda)=\{17+9\omega, 23+14 \omega\}$.
To determine the correct edge orientation,  we first determine the color of the edge using
the first component of the midpoint ($(17+23)/2=20$) and the first component of the difference ($(23-17)=6$) to obtain
a color $c(\{17+9\omega,23+14\omega\})=\lfloor20/6\rfloor\,(mod\, 2) = 1$. 
Since the color is $1$, the closer endpoint becomes the channel 2 label. Hence $18+10\omega$ is assigned
the label $(23+14\omega,17+9\omega)$.

We now illustrate the decoding procedure. Assume that
$\oer=(23+14\omega,17+9\omega)$ and write its undirected
version using the basis vectors of the sublattice $\bu$ and
$\bv$ to get $\{4\bu+3\bv,3\bu+2\bv\}$. Look for an equivalent
edge in Table~\ref{table-alphau}. One such edge is
$\{O,B\}=\{0,\bu+\bv\}$. Determine the shift required to make the
edges coincide. In this case
$\{4\bu+3\bv\}-\lambda'=\{\bu+\bv,0\}$, with
$\lambda'=3\bu+2\bv$.
Upon looking up Table~\ref{table-alphau}, we
find that lattice point $b$ with representation $1+\omega$ has the undirected edge $\{O,B\}$ as label.
We shift this point by adding $\lambda'=3\bu+2\bv$ to get one of the candidate
points $\lambda=3\bu+2\bv +(1+\omega)=18+10\omega$. The
other candidate point is obtained from the Property 3 of the
labeling function (the sum of the endpoints of an edge is equal to the sum of
the points that it labels), and is $22+13\omega$. Observe that there is
another edge in Table~\ref{table-alphau} which is equivalent to
the edge we wish to decode. We would obtain exactly the same set
of
candidate points if we used this edge. In order to determine the
correct point, since $c(e)=1$, we decode to  the point which is closer to the channel 2 label, namely $18+10\omega$.

Several observations can be made at this point. In
Table~\ref{table-alphau} there are two kinds of undirected edges (of positive
length)--those which are diameters of a circle centered at $0$ and those which are not.  
The diameters are the edges $\{A,D\}$, $\{B,E\}$ and $\{C,F\}$. Both orientations of a diameter are used  to
label points in $V_0(0)$, whereas only one orientation of a non-diameter is used.
For an edge which is not a diameter, the remaining orientation labels
a lattice point outside $V_0(0)$ as determined by the shift property. For example, consider the label $\{L,C\}$.
The directed label $(L,C)$ is the label for the point $ac=1-2\omega \in V_0(0)$. The lattice point which is
labeled by $(C,L)$ is given by $2\mu -(1-2\omega)=5-\omega-(1-2\omega)=4+\omega$, which belongs to the discrete
Voronoi set of the sublattice point $A$.

Notice that the labeling function shown in
Table~\ref{table-alphau} exhibits an additional symmetry that we
have so far not used. If a point is rotated by a multiple
of $\pi/3$ radians about the origin, its corresponding undirected
label is also rotated by the same amount about the origin.
Consider $\Gamma=\{\gamma_i=\exp{ik\pi/3},k=0,1,\ldots,5\}$, a
rotation group of order 6. By considering equivalence classes
relative to this group, we can reduce the size of the
table by listing only the undirected edges for the points
$\{o,a,g,n,t,s\}$.

The above example illustrates the basic steps that are to be followed in order to
label the points in $\Lambda$. Additional details and some
underlying theory is presented next.

\subsection{\bf General Principles}

The construction of our labeling function involves the following steps.
\begin{enumerate}
\item Selection of a geometrically similar sublattice of given index $N$.
\item Construction of $V_0(0)$, the discrete Voronoi set around $0$.
\item Establishing a mapping between elements of $V_0(0)$ and undirected
edges in such a way that certain constraints are
satisfied. The optimal construction requires that a specific linear programming problem be solved.
\item Extension of the mapping to the entire lattice.
\item Identification of a specific directed edge to associate with a lattice point,
once the undirected edge is known.
\end{enumerate}
In the remainder of this section, we will further describe items 3, 4 and 5. But first we state 
the following guiding principle.

\subsubsection{A Guiding Principle}
\begin{figure}[htb]
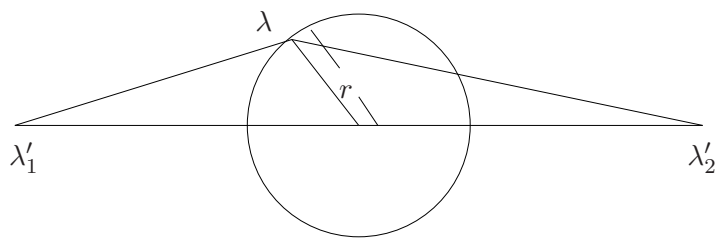

\begin{center}
\input simple_config.pstex_t
\end{center}
\caption{A lattice point $\lambda$ and its label $(\lambda'_1,\lambda'_2)$. $\lambda$ is at distance $r$
from the midpoint of the labels.}
\label{fig-simple_config}
\end{figure}
Suppose that $\lambda$ receives the label $\oer=(\lambda'_1,\lambda'_2)$ as illustrated in Fig.~\ref{fig-simple_config}. 
Then $d_s(\overrightarrow{e})$ satisfies the identity
\begin{equation}
2d_s(\oer)=\|\lambda-\lambda'_1\|^2+\|\lambda-\lambda'_2\|^2=(1/2)\|\lambda'_1-\lambda'_2\|^2+2r^2,
\label{eqn-a2-2}
\end{equation}
where $r^2=\|\lambda-(\lambda'_1+\lambda'_2)/2 \|^2$.
This identity is known as the parallelogram
law (see for example \cite{Auslander:1}, page 3).
From this we infer that in order to keep $\bar{d}_s(\overrightarrow{e})$ small, $\lambda'_1$ and $\lambda'_2$ should be as close together
as possible, and $\lambda$ should be as close as possible to their
midpoint. This leads to the following:
\par\bigskip
\fbox{\parbox{6.0in}{%
{\bf Guiding Principle}: Choose the shortest possible edge with midpoint as close as possible to the point to be
labeled.}}

\subsubsection{\bf Optimal Construction of $\calE_u(0)$ and $\alpha_u$ (Items 3 and 4)}
\label{sec-alphau}

\begin{figure}[htb]
\centerline{\psfig{figure=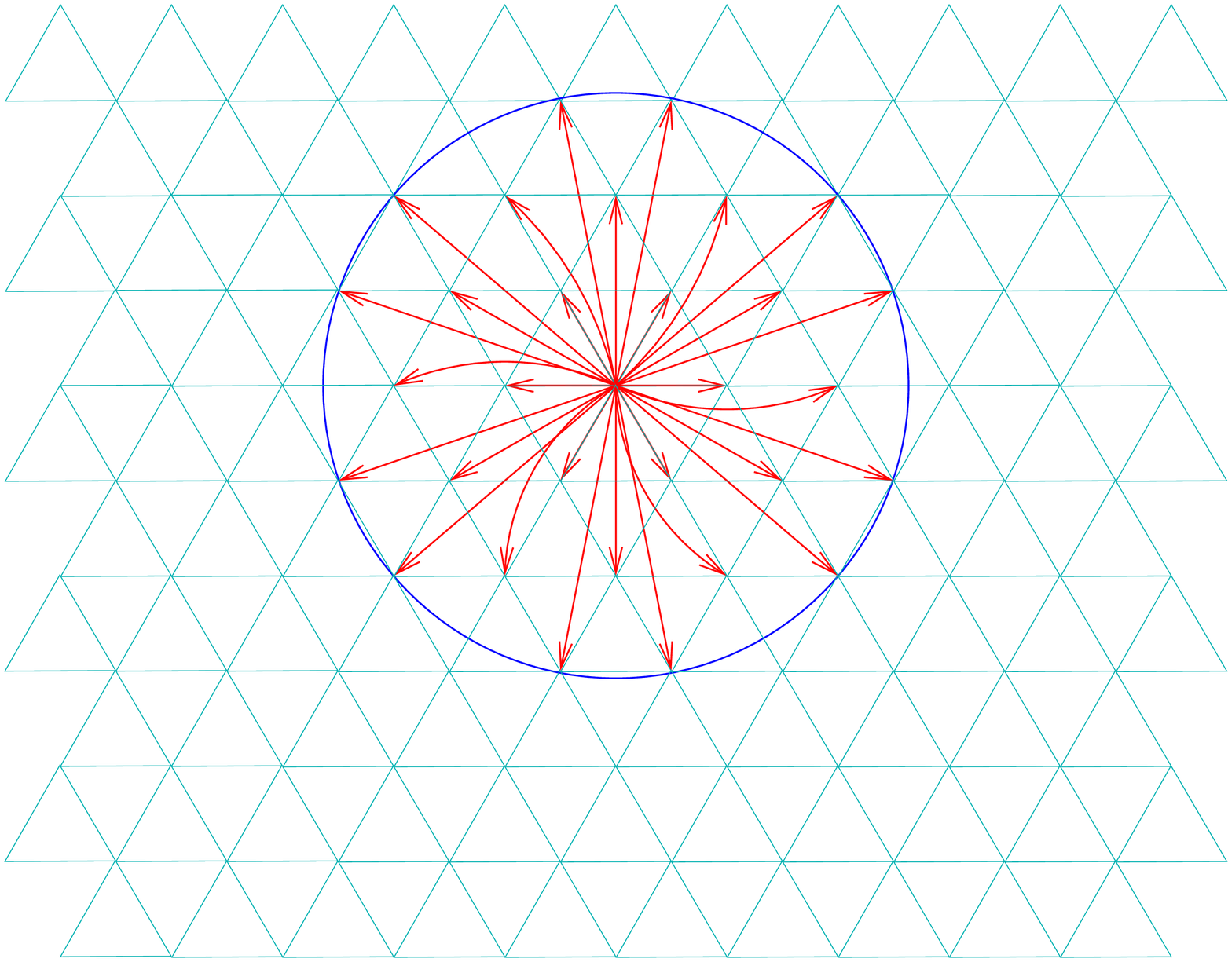,width=3in}}
\caption{The set $\calE_d|_{\alpha_1=0}$ for $N=31$.}
\label{fig-e0star}
\end{figure}

\begin{figure}[htb]
\begin{center}
\input{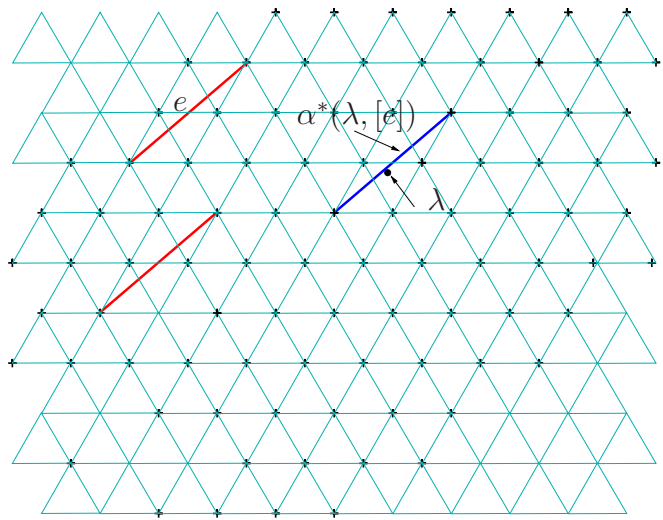}
\end{center}
\caption{Illustrating the optimal edge selection within a given equivalence class. Shown is an edge $e$, another
edge in its equivalence class $[e]$ and the optimal edge in this class for labeling the point $\lambda$. The points
on the hexagonal grid are sublattice points.}
\label{fig-optimal-edge}
\end{figure}

The starting point for determining the set of edges is to compute
$\calE_d|_{\alpha_1=0}$, which is chosen to be the set of $N$ shortest undirected sublattice edges of the form
$(0,\lambda')$. This set is illustrated in
Fig.~\ref{fig-e0star}. The corresponding set $\calE_u|_{\alpha_1=0}$ is
obtained by replacing each directed edge with its undirected version.
Observe that
edges of positive length in $\calE_u|_{\alpha_1=0}$ occur in equivalent pairs.
Next, for
every $\lambda \in V_0(0)$ and $e \in
\calE_u|_{\alpha_1=0}$, determine $d_s(\lambda,[e]) = \min_{\tilde{e}\in [e]}d_s(\lambda,\tilde{e})$ and let
$\alpha^*(\lambda,[e])$ be an edge $\tilde{e} \in [e]$ which achieves
this minimum.
In other words,
$\alpha^*(\lambda,[e])$ is the ``closest'' edge to $\lambda$ which lies in
the same equivalence class as $e$. Clearly, from the guiding principle, the ``closest''
edge will be the one whose midpoint is closest to $\lambda$.
The edge $\alpha^*(\lambda,[e])$ for a given $\lambda$ and
$e$ is illustrated in Fig.~\ref{fig-optimal-edge}.
\clearpage

Consider all one-to-one maps
$\beta~:~V_0(0)\rightarrow
\calE_u|_{\alpha_1=0}$
which satisfy the
constraint that equivalent points are mapped to equivalent edges.
From among all such maps $\beta$, choose
$\beta^*$ so as to minimize $\sum_{\lambda \in V_0(0)}d_s(\lambda,[\beta(\lambda)])$.
The map $\beta^*$ sends a lattice point to an edge coset
element in $\calE_u|_{\alpha_1=0}$ in an optimal way, thus identifying the
best edge coset for a given lattice point.
Let $\alpha_u(\lambda)=\alpha^*(\lambda,[\beta^*(\lambda)])$.
Since $\beta^*$ identifies the optimal edge
coset for each lattice point and $\alpha^*$ identifies the best
coset representative, given the lattice point and the coset, we
obtain the optimal edge for each lattice point by composing these
two mappings.
It follows that $\calE_u(0)=\{\alpha_u(\lambda)~:~\lambda \in V_0(0)\}$.
We extend the mapping to the lattice $\Lambda$ using the rule
$\alpha_u(\lambda+\lambda')=\alpha_u(\lambda)+\lambda'$.

\begin{figure}[htb]
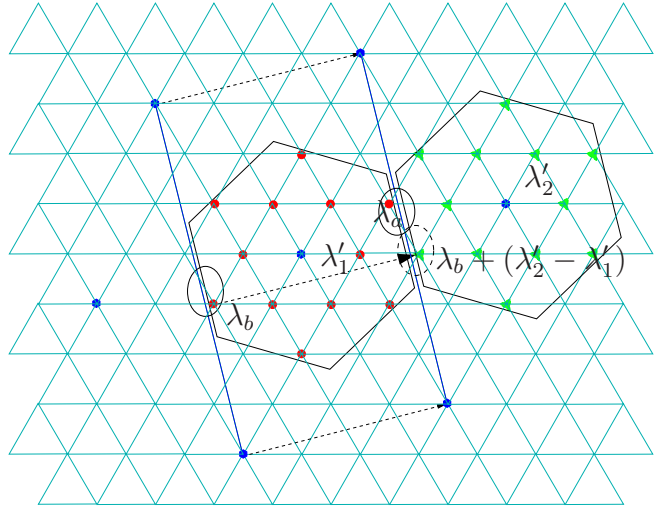

\begin{center}
\input ToExplainConstraint.pstex_t
\end{center}
\caption{Explanation of reason for mapping equivalent points in $V_0(0)$ to equivalent edges. The lattice $A_2$ together with
a similar sublattice of index $N=13$ is shown. Points $\lambda_a$ and $\lambda_b$ in $V_0(\lambda_1')$ are equivalent and
$\lambda_a+\lambda_b=2\lambda_1'$. $\lambda_a$ and $\lambda_b$
are mapped to equivalent edges $e_a$ and $e_b$. Here $e_a=e_b+(\lambda_2'-\lambda_1')$. The point $\lambda_b+(\lambda_2'-\lambda_1')$ lies in $V_0(\lambda_2')$ and is
labeled by $e_b+(\lambda_2'-\lambda_1')=e_a$. Thus $\lambda_a
+(\lambda_b+(\lambda_2'-\lambda_1'))=\lambda_1'+\lambda_2'$. The midpoint of two points labeled by the same undirected edge
is equal to the midpoint of the edge. }
\label{fig-explain-constraint}
\end{figure}

The constraint imposed on the mapping $\beta$ needs some explanation.
It arises from the third assumption that we made about the labeling function,
namely, that the sum of the end points of an edge is equal to the sum of the two
points that it labels. By requiring that equivalent points map to equivalent edges,
it can be shown that the midpoint of two lattice points that share the same (undirected) label
coincides with the midpoint of the label itself. A graphical justification for this is provided in
Fig.~\ref{fig-explain-constraint}.
More formally, the argument is as follows. The points in $V_0(0)$ occur in
equivalent pairs. If two points in $V_0(0)$ are equivalent, they sum to $0$. Consider
the pair $\lambda_a$ and $-\lambda_a$ and an edge $e$. If the
edge in $[e]$ which is closest to $\lambda_a$, say $e_a$, has midpoint $a$, then the edge in $[e]$ which is closest
to $-\lambda_a$ will have midpoint $-a$. Thus $\calE_u(0)$ contains the
edges $e_a$ and $-e_a$ (which may be identical).
Now from the shift property, $e_a$ also lies
in $\calE_u(2a)=\calE_u(0)+2a$ (note that $2a$ is a sublattice point) and the point it labels,
say $\lambda_b$, is given by $\lambda_b=-\lambda_a+2a$. Thus the two points that are labeled using
the undirected edge $e_a$, namely $\lambda_a$ and $\lambda_b$, satisfy $\lambda_a +
\lambda_b=2a$, i.e., the midpoint of the edge coincides with the midpoint of the points that
receive this edge as label. To summarize, the constraint is a sufficient condition
to ensure that the labeling function has Property 3. 

It is to be noted that the optimal mapping $\beta^*$ can be obtained using
standard techniques from linear programming~\cite{LUENI}. Also
observe that if we define the group $\Gamma=\{1,-1\}$, 
we can force $\beta^\ast$ to satisfy the constraint
by considering
only $V_0(0)/\Gamma$ and $\calE_u|_{\alpha_1=0} / \Gamma$. This
problem is one of matching $(N-1)/2$ lattice points to $(N-1)/2$
edge classes. Further reductions in complexity may be obtained by
using a larger group $\Gamma$, in which case the problem is 
reduced to matching $(N-1)/o(\Gamma)$ points to $(N-1)/o(\Gamma)$
edge classes,
where $o(\Gamma)$ is the order of the group.
In the case of the $A_2$ lattice, for example, we could take $\Gamma$
to be a group of order 6.

The mapping $\alpha_u$ identifies the edge to be associated with $\lambda$
up to the orientation of the edge. The correct orientation is
determined by an edge orientation rule, designed to
maintain balance between the two descriptions.
This is described next.

\clearpage
\subsubsection{\bf Balance}
\label{sec-balance}
\begin{figure}[htb]
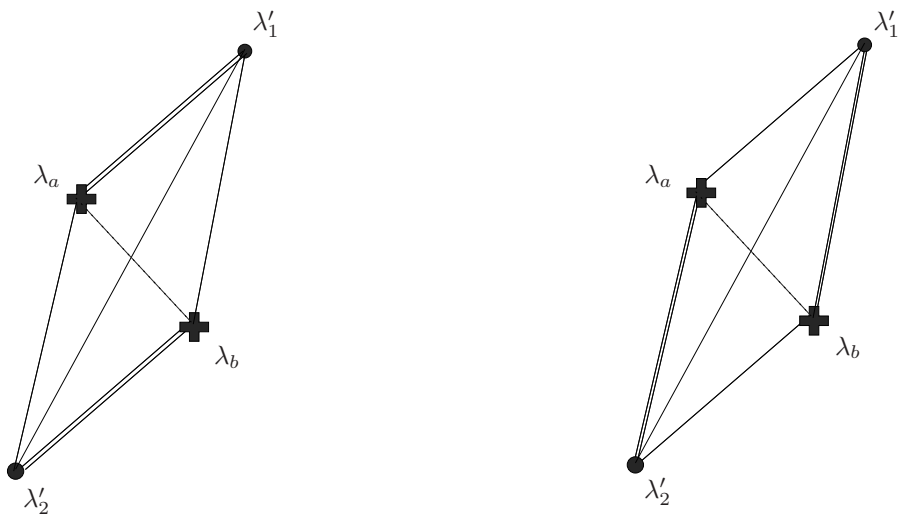

\begin{center}
\input balance.pstex_t
\end{center}
\caption{
Given an edge $e = \{ \la'_1, \la'_2 \}$ and the two lattice points $\{ \la_a, \la_b \}$ that it labels, with the midpoint of $e$ coinciding with the midpoint
$(\la_a + \la_b )/2$, there are two possibilities for the selection rule.
In the figure the lattice points are joined to the sublattice points by both a single line,
indicating the first component, and a double line, indicating the second component.
On the left the double lines are shorter than the single lines, a selection rule which favors the second description, that is, $d_2$ is smaller than $d_1$.
On the right the single lines are shorter than the double lines, a rule which favors the first description.
Balance is achieved by requiring that both rules are used equally often.}
\label{fig-balance}
\end{figure}

\begin{figure}[htb]
\centerline{\psfig{figure=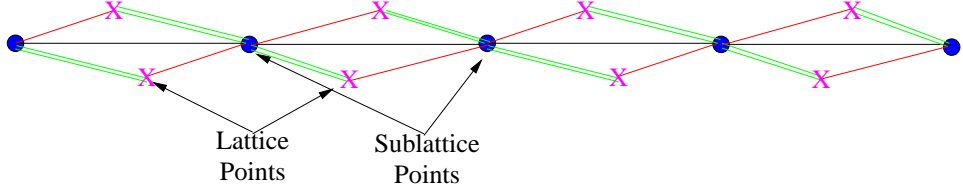,width=5in}}
\caption{An illustration of the alternating structure along
a straight line for maintaining balance between the two
descriptions. A single line connecting a lattice point to a sublattice point indicates the channel 1 label for that lattice point, a
double line indicates the channel 2 label.}
\label{fig-alternating}
\end{figure}

Given an edge $e$ and two points $\{\lambda_a, \lambda_b \}$ that it labels,
there are two ways to establish a one-to-one
correspondence between the points and the two directed edges $\{\oer,\oel\}$,
as described in Fig.~\ref{fig-balance}.
The first selection rule favors the second
description, the second favors the first description.
Note that the
distortions are anti-symmetric, i.e., $d_1(e)$ with the first selection rule is equal to $d_2(e)$
for the second, and $d_2(e)$ for the first rule is equal to $d_1(e)$ for the second. Because of this
anti-symmetry, balance is
attainable in an average sense by ensuring that the two correspondences are used equally often.
We achieve this by requiring that the two correspondences are alternated along any straight line of the
edge graph, as shown in Fig.~\ref{fig-alternating}.
This is the purpose of the edge coloring rule which is described next.

\subsubsection{\bf An Edge  Coloring Rule}
\label{sec-coloring}
Each edge $e$ is assigned a bit (or color\footnote{During this research we used red and green to color the edges.}) $c(e)$ as follows.
Consider an edge of the form $(a+b\omega,c+d\omega)$. Let $\Delta_1=|c-a|$ and let $\Delta_2=|d-b|$.
Then 
\begin{equation}
c(e)=\left\{\begin{array}{ll}
0, & 
\mbox{if $\Delta_1>0$ and  $\lfloor(c+a)/(2\Delta_1) \rfloor$ is even},\\
 & \mbox{or if $\Delta_1=0$ and $\lfloor(d+b)/(2\Delta_2) \rfloor$ is even},\\
1, & \mbox{otherwise}.
\end{array}\right.
\label{eqn-ColoringRule}
\end{equation}
This coloring rule ensures that adjacent edges along any
straight line have a different color. Note that the color does not depend on the orientation of
an edge, i.e., $\oer$, $\oel$ and $e$ have the same color. 
The coloring rule is an ingredient in the edge and point selection rules that we now define.

\subsubsection{\bf Edge Direction Rule and Point Selection Rule}
\label{section-edge-direction}
Given an (undirected) edge $e$ and a point $\lambda$ for which this edge is a label, we choose
an orientation or direction for the edge using a rule that depends on the color of the edge.
Let $e=\{\lambda'_1,\lambda'_2\}$ and let $\mu=(\lambda'_1+\lambda'_2)/2$.
The two  rules $s_c(e,\lambda)$, where $c=0$ or 1 is the color of the edge, are defined as follows
($\times$ denotes the cross- or vector-product of two vectors and $\langle
\cdot,\cdot\rangle$ their inner product):
\begin{equation}
s_0(e,\lambda)=\left\{\begin{array}{ll}
(\lambda'_1,\lambda'_2), & \mbox{$\langle\lambda'_1-\lambda'_2,\lambda-\mu\rangle>0$ or}\\
&  \mbox{$\langle\lambda'_1-\lambda'_2,\lambda-\mu\rangle=0$ and $\mbox{sign}((\lambda'_1-\lambda'_2)\times (\lambda-\mu))>0$,}\\
(\lambda'_2,\lambda'_1), & \mbox{otherwise}.
\end{array}
\label{eqn-DirectionRuleA}
\right.
\end{equation}
and
\begin{equation}
s_1(e,\lambda)=\left\{\begin{array}{ll}
(\lambda'_2,\lambda'_1) & \mbox{$\langle\lambda'_1-\lambda'_2,\lambda-\mu\rangle>0$ or}\\
&  \mbox{$\langle\lambda'_1-\lambda'_2,\lambda-\mu\rangle=0$ and $\mbox{sign}((\lambda'_1-\lambda'_2)\times (\lambda-\mu))>0$,}\\
(\lambda'_1,\lambda'_2) & \mbox{otherwise}.
\end{array}
\right.
\label{eqn-DirectionRuleB}
\end{equation}
Observe that the result of either rule is the same whether we write $e=\{\lambda'_1,\lambda'_2\}$ or $e=\{\lambda'_2,\lambda'_1\}$.

For decoding, since two lattice points receive a  label from a given undirected edge, we need to be able to tell which
point is being labeled, given the edge orientation. This is the reverse of the edge direction rule. 
Thus, given a directed edge $\oer=(\lambda'_1,\lambda'_2)$ with midpoint $\mu$ and a lattice point 
$\lambda$ which could have received this label, the {\em Point Selection Rule} $g_c(\oer,\lambda)$ selects $\lambda$ or $2\mu-\lambda$ 
based on the edge color $c$, and is given by
\begin{equation} 
g_0(\oer,\lambda)=\left\{ \begin{array}{ll}
\lambda & \mbox{$\langle \lambda'_1-\lambda'_2,\lambda-\mu \rangle>0$ or}  \\
          & \mbox{$\langle \lambda'_1-\lambda'_2,\lambda-\mu \rangle=0$ and $\mbox{sign}((\lambda'_1-\lambda'_2)\times (\lambda-\mu))>0$,}  \\
2\mu-\lambda & \mbox{otherwise}.
\end{array}
\right.
\label{eqn-PointSelRuleA}
\end{equation}
and
\begin{equation} 
g_1(\oer,\lambda)=\left\{ \begin{array}{ll}
2\mu-\lambda & \mbox{$\langle \lambda'_1-\lambda'_2,\lambda-\mu \rangle>0$ or}  \\
          & \mbox{$\langle \lambda'_1-\lambda'_2,\lambda-\mu \rangle=0$ and $\mbox{sign}((\lambda'_1-\lambda'_2)\times (\lambda-\mu))>0$,}  \\
\lambda & \mbox{otherwise}.
\end{array}
\right.
\label{eqn-PointSelRuleB}
\end{equation}

\subsubsection{\bf Constructing the map $\alpha$}
Using the notation previously established, we obtain $\alpha$ as follows.
Given $\lambda$, let $e=\alpha_u(\lambda)$ and let $c=c(e)$. Then $\alpha(\lambda)=
s_c(e,\lambda)$.

\subsubsection{\bf Proof that the Reuse Index is Correct}
We need to show that for any sublattice point $\lambda'$,
exactly $N$ lattice points have a label of the form
$(\lambda',*)$ and exactly $N$ lattice points
have a label of the form $(*,\lambda')$.

Consider the graph drawn on the vertex set $\Lambda'$, where the edge set is
$\calE_u=\bigcup_{\lambda' \in \Lambda'}(\calE_{u|\alpha_1=0}+\lambda')$. This is
the entire set of undirected edges or labels, or equivalently, the range of the map $\alpha_u$.
Observe that for any $\lambda' \in \Lambda'$, there are exactly $N$ edges in $\calE_u$
of the form $\{\lambda',*\}$, or equivalently, there are exactly $N$ edges incident
on the vertex $\lambda'$. Each edge $\{\lambda',\lambda''\}$ of positive length labels two points,
of which one
receives the label $(\lambda',\lambda'')$, the other $(\lambda'',\lambda')$. Thus
exactly $N$ lattice points receive a label of the form $(\lambda',*)$ and exactly
$N$ receive a label of the form $(*,\lambda')$.

\subsubsection{\bf Further Reduction in Complexity: Group Construction}
Instead of imposing the constraint on the map $\beta$
that was used in Section 4.2.2, we could alternatively regard this
map
as an {\em unconstrained} map between cosets $\Lambda/\Lambda'/\Gamma$
and $\calE_u / \Lambda' / \Gamma$, where $\Gamma$ is the group of rotations $\{1,-1\}$ (in
complex notation since we are talking about $A_2$).
Now we need only establish
a correspondence between two sets of size $(N-1)/2$, where $2$ is the order of the group.
Further reductions in complexity arise from selecting a larger group that contains the
group $\{1,-1\}$ as a subgroup.
For $A_2$, we used the group  of rotations
$\{\exp(ik\pi/6),k=0,1,2,3,4,5\}$. This reduces the complexity of matching problem
to sets of size $(N-1)/6$. Precise conditions that the group must satisfy, and 
further reasons for using a group, are
explained in Section 6.

\subsubsection{Numerical Results}
\begin{figure}[htb]
\centerline{\psfig{figure=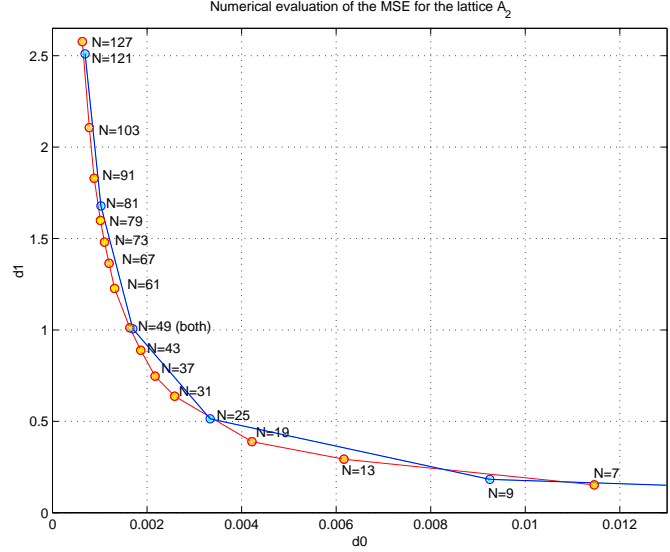,width=4in}}
\caption{A plot of $\bar{d}_1$ v.s. $\bar{d}_0$ for the hexagonal
lattice $A_2$ and the integer lattice $\ZZ$.}
\label{fig-A2-results}
\end{figure}
The results of our optimization procedure
for the hexagonal lattice $A_2$ are
displayed in Fig.~\ref{fig-A2-results}, along with comparisons with
the $\ZZ$ lattice. These results have been obtained for a uniformly
distributed memoryless source by computing the optimal labeling function and then
evaluating the expression 
(obtained from (\ref{eqn-app-5}) for a uniform probability density function)
\begin{equation}
\dbar_s=\dbar_0 +(1/N)\sum_{\lambda \in V_0(0)}d_s(e),
\label{eqn27}
\end{equation}
where $e$ is the edge that labels $\lambda$, $\bar{d}_0$ is computed
using (\ref{eqn-d0}) and known values of $G(\Lambda)$ \cite{SPLAG}.
To be comparable with $A_2$, the value of $N$ shown for
the $\ZZ$ lattice is the square of the actual reuse index for dimension $L=1$. 
To be more specific, let $N_{A_2}$ denote the index for the $A_2$ sublattice, and let
$N_{\ZZ}$ be the index of the sublattice of $\ZZ$. 
Then the value of $N$ stated in the figure is $N_{\ZZ}^2$ for $\ZZ$ and $N_{A_2}$ for $A_2$. 
Also it should be noted that for each $N$, both lattices have been scaled in
order to keep the product $N\nu$ constant, where $\nu$ is the volume of a
fundamental region of the {\em scaled} lattice. From (\ref{eqn-rate-1}), this is equivalent to 
keeping the rate constant for all points on the graph.
It is seen that small performance improvements
are obtained by using $A_2$ instead of $\ZZ$. Another benefit is that we obtain many more points
in a given interval of (say) the side distortion, compared with 
the $\ZZ$ lattice. 

\section{Labeling Functions for General Lattices} 
\label{sec-gen}
In this section we describe how to label a general lattice $\Lambda$ using a sublattice
$\Lambda'$ of index $[\Lambda:\Lambda']=N$. The basic steps remain the same as for $A_2$.
The main differences arise in the selection of the sublattice and in the use of a group to simplify the construction.
We begin by establishing certain general conditions that the group should satisfy.
Specific groups and sublattices will then be given for certain particular
lattices.
We will use $\bG$ to denote a generator matrix for $\La$ and $c \tilde{\bG} \bG$
for the generator matrix for the similar sublattice $\La'$, where $c$ is an
scalar and $\tilde{\bG}$ is a unitary matrix.
Our convention is that the columns of a generator matrix are a basis for the lattice.

\subsection{Sufficient Conditions for the Group, $\Gamma$}
The smallest group we can use is $\Gamma=\{I_L,-I_L\}$, where $I_L$ is the
$L$-dimensional identity matrix.   Here we show why a group is
useful for reducing the size of the optimization problem and
derive certain conditions that the group should satisfy.

The motivation for using a group is to make use of
inherent symmetries in the lattice and sublattice.
Our objective is to
partition the discrete Voronoi set and the edge set into subsets 
of equal size with certain distance properties. More specifically, for any
subset of lattice 
points $\{\lambda_i,~i=1,2,\ldots,M\}$ in this partition and any
subset of edges $\{e_i,~i=1,2,\ldots,M\}$ in the partition the
set of distances ${\cal D}_i=\{d_s(\lambda_i,e_j),~j=1,2,\ldots,M\}$
should be independent of $i$. Such sets of points and edges can be
obtained by identifying a group of transformations
$\Gamma=\{\gamma_k,~k=1,2,\ldots,M\}$ and then taking the 
members of the partition to be orbits under the action of this
group. The group we look for should have following properties.
\begin{enumerate}
\item $\Gamma$ contains $-I_L$.
\item $\Gamma$ is an orthogonal group.
\item $\Gamma$ preserves the lattice $\Lambda$.
\item  $\Gamma$ acts fixed-point free on $\Reals^L$, i.e., for
any $\gamma \in \Gamma$, $\gamma$ not the identity, $\gamma x =x
\implies x=0$.
\item The order of the group divides the g.c.d. of all the shell
sizes.
\item $\Gamma$ preserves the sublattice $\La'$.
\end{enumerate}
Property 2 ensures that orbits of the group lie entirely within a
shell of the lattice. Property 3 is obvious, otherwise orbits would
contain points other than lattice points. Property 4 ensures that
all orbits are of equal size. Property 5 ensures that each shell
is partitioned into an integral number of orbits. Property 6
is similar to Property 3.

We now look at Properties 3 and 6 more closely. Every lattice point
$\lambda$ can be written in terms of the generator matrix $\bG$
and an integer vector $u$ as $\lambda=\bG u$. By
requiring that $\gamma_i\bG u=\bG\gamma_ju$ for some $\gamma_j \in \Gamma$,
we ensure that the lattice is preserved. Similarly, if the 
sublattice has generator matrix $c\tilde{\bG}\bG$, where $c$ is a scalar,
then
by requiring $\gamma_i\tilde{\bG}=\tilde{\bG}\gamma_j$, for some
$\gamma_j \in \Gamma$, the sublattice will also be preserved.
In other words we require that $\Gamma$ {\em normalizes} $\bG$ and
$\tilde{\bG}$.
\subsection{Group Construction and Generator Matrices}
For the space lattice $\ZZ^2$, we take $N$ to be an odd number of the form
$a^2 + b^2$ where $a,b \in \ZZ$.
(It is shown in \cite{ConRaiSlo1} that similar sublattices of $\ZZ^2$ exist whenever $N$ is a sum of two squares.)
The generator matrices are $\bG=I_2$ and
$\bG'=\tilde{\bG} \bG$, where
\begin{equation}
\tilde{\bG}=\left(
\begin{array}{cc}
a & -b \\
b & a
\end{array}
\right).
\end{equation}
For the group we use
\begin{equation}
\Gamma=\left\{\pm I_2,
\pm\left(
\begin{array}{cc}
0 & -1 \\
1 & 0
\end{array}
\right)
\right\},
\end{equation}
a group of order 4.

For the lattice $\ZZ^4$
we take $N$ to be an odd perfect square (again see \cite{ConRaiSlo1}).
Let
$N=a^2+b^2+c^2+d^2$, $a,b,c,d \in \ZZ$ (any integer can be written this way).
The generator matrices are $I_4$ and $\bG'=\tilde{\bG} \bG$, with
\begin{equation}
\tilde{\bG}=\left(
\begin{array}{cccc}
a & -b & -c & -d \\
b & a & -d & c \\
c & d & a & -b \\
d & -c & b & a
\end{array}
\right).
\end{equation}
The group $\Gamma$ is
\begin{equation}
\Gamma=\left\{\pm I_4,
\pm\left(
\begin{array}{cccc}
0 & -1 & 0 & 0 \\
1 & 0  & 0 & 0 \\
0 & 0 & 0 & 1 \\
0 & 0 & -1 & 0
\end{array}
\right),
\pm\left(
\begin{array}{cccc}
0 & 0 & -1 & 0 \\
0  & 0  & 0 & -1 \\
1 & 0 & 0 & 0 \\
0 & 1 & 0 & 0
\end{array}
\right),
\pm\left(
\begin{array}{cccc}
0 & 0 & 0 & -1 \\
0  & 0  & 1 & 0 \\
0 & -1 & 0 & 0 \\
1 & 0 & 0 & 0
\end{array}
\right)
\right\},
\end{equation}
a group of order 8.

For $\ZZ^8$, it is easiest to start by specifying the group.
Let
\begin{equation}
\gamma_1=\left(
\begin{array}{rrrrrrrr}
0 & 1 & 0 & 0 &  & & &\\
0  & 0  & 1 & 0 &  & & &\\
0 & 0 & 0 & 1 &  & & &\\
-1 & 0 & 0 & 0 &  & & &\\
 & & & & 0 & 1 & 0 & 0\\
 & & & & 0 & 0 & 1 & 0 \\
 & & & & 0 & 0 & 0 & 1\\
 & & & & -1 & 0 & 0 & 0
\end{array}
\right)
\end{equation}
and
\begin{equation}
\gamma_8=\left(
\begin{array}{rrrrrrrr}
 & & & & 1 & 0 & 0 & 0  \\
& & & & 0 & 0 & 0 & -1 \\
& & & & 0 & 0 & -1 & 0 \\
& & & & 0 & -1 & 0 & 0 \\
-1 & 0 & 0 & 0 & & & & \\
0 & 0 & 0 & 1 & & & & \\
0 & 0 & 1 & 0 & & & & \\
0 & 1 & 0 & 0 & & & & 
\end{array}
\right).
\end{equation}
Then we take the group to be
$\Gamma=\{\gamma_1^i,\gamma_8 \gamma_1^i,~i=0,1,2,\ldots,7 \}$,
a group of order 16.
The generators
for the lattices are $\bG=I_8$ and $\bG' = \tilde{\bG} \bG$, where the
$i$th column of $\tilde{\bG}$ is $\gamma_i v$,
where $v= (a~0~b~0~c~0~d~0)^{tr}$, $a^2 + b^2 + c^2 + d^2 =N$, the index of
$\La'$, and $i=0,1, \ldots, 7$.

\begin{figure}[htb]
\centerline{\psfig{figure=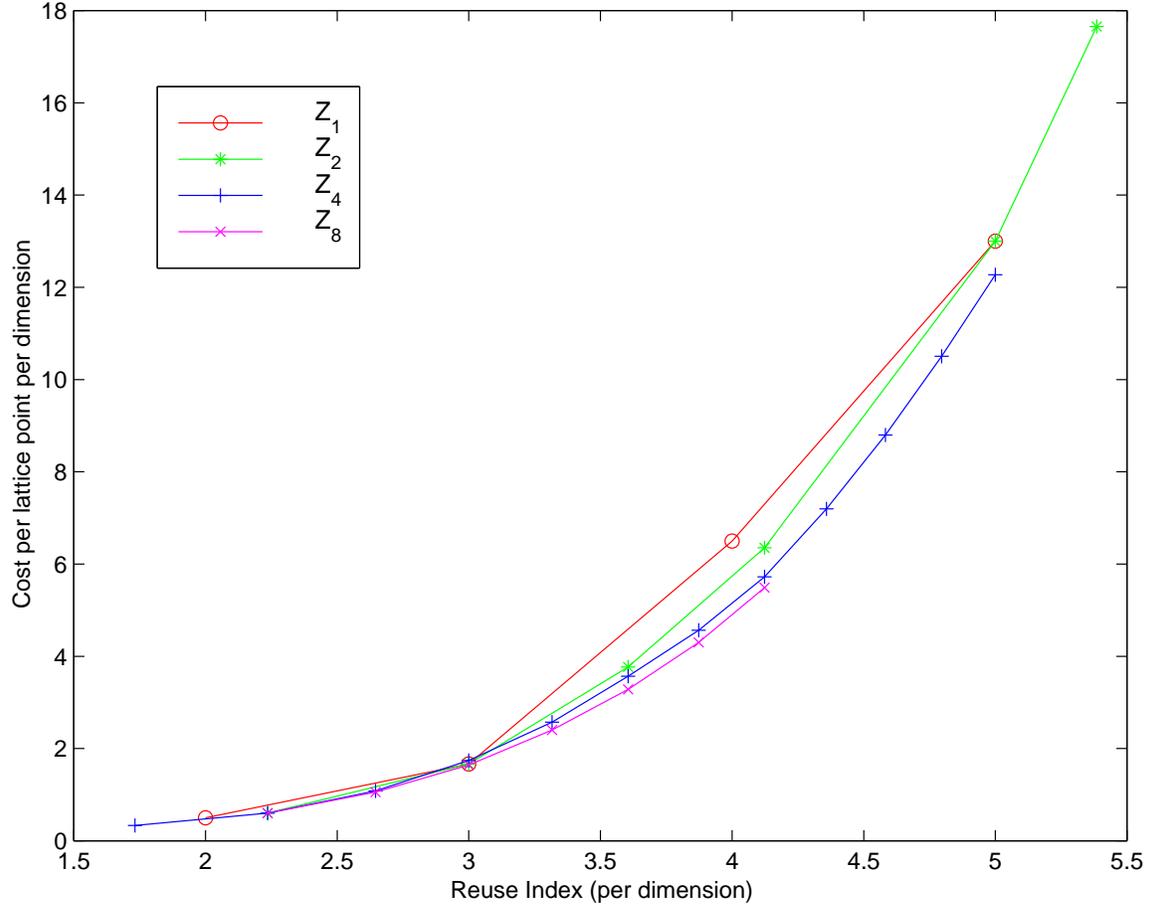, width=6in}}
\caption{A plot of $\bar{d}_s$ v.s. the reuse index (per
dimension) $N^{1/L}$ for
the lattices $\ZZ^i$, $i=1,2,4,8$.}
\label{plot-Zi}
\end{figure}
Numerical computations for the lattices $\La = \ZZ^i$
(standardized to have minimal length 1) for $i=1,2,4,8$,
are presented in Fig.~\ref{plot-Zi}. These results have also been obtained
for a uniformly distributed memoryless source. Since the two-channel
distortions are identical for all the lattices considered in this figure, we
have only plotted the excess distortion term $(1/N)\sum_{\lambda \in V_0(0)}d_s(\lambda)$
against the per-dimension reuse index $N^{1/L}$.

\section{\bf Asymptotic Analysis}
\label{sec-asymp}
We now derive upper and lower bounds on the distortion $\bar{d}_s$ as given
by (\ref{eqn-app-5}).
As we have already seen in (\ref{eqn27}),
the regularity of the labeling function and the  high rate
assumption lead to the following simplification
for the the expression in (\ref{eqn-app-5}): 
\begin{equation}
\bar{d}_s=\bar{d}_0 + \frac{1}{N}\sum_{\lambda \in V_0(0)}d_s(e),
\label{eqn-a2-1}
\end{equation}
where $e$ is the label for the lattice point $\lambda$.
Thus in order to analyze $\bar{d}_s$ we only need
to consider edges that label points
in the discrete Voronoi set $V_0(0)$. Our analysis relies on
a precise knowledge of the lengths of these
edges.
Suppose that $\lambda$ is labeled by
the edge $e=(\lambda'_1,\lambda'_2)$,
as in Fig.~\ref{fig-simple_config}. 
Let $l(e)$ denote the length of the edge
$e=\{\lambda'_1,\lambda'_2\}$, i.e.,
$l(e)=\|\lambda'_1-\lambda'_2\|$, and let $r$ be the distance of
the point $\lambda$
from the midpoint $(\lambda'_1+\lambda'_2)/2$. Then 
\begin{equation}
2d_s(e)=\|\lambda-\lambda'_1\|^2+\|\lambda-\lambda'_2\|^2=(1/2)l(e)^2+2r^2.
\label{eqn-a2-2prime}
\end{equation}
A simple lower bound for $d_s(e)$ is obtained by setting $r$ to
zero in the above equation. In order to obtain an upper bound,
observe that the midpoint of a sublattice edge can always be made to lie
in $V_0(0)$ by a suitable sublattice shift.  An
upper bound is then obtained
(using the triangle inequality) by replacing $r$ with
$r^*(\Lambda')=2R(\Lambda')$, where $R(\Lambda')$ is the
covering radius of the sublattice.
Thus we have
\begin{equation}
(1/4)l^2(e)\leq d_s(e) \leq (1/4) l^2(e)+ r^*(\Lambda')^2,
\end{equation}
which by using (\ref{eqn-a2-1}) leads to the bounds
\begin{equation}
\bar{d}_0 + \frac{1}{4N} \sum_{\lambda \in V_0(0)}l^2(e) \leq
\bar{d}_s \leq \bar{d}_0
+ \frac{1}{4N} \sum_{\lambda \in V_0(0)}l^2(e) + r^*(\Lambda')^2.
\label{eqn-a2-5}
\end{equation}

\subsection{Asymptotic Performance}
\label{sec-asymp-gen}
In order to carry out a rate-distortion analysis, it is necessary to scale $\Lambda$ and
$\Lambda'$ by a real number $\beta$. We will
use $\nu(\beta \Lambda)$ to denote the volume of a fundamental region for the scaled lattice $\beta
\Lambda$. Clearly $\nu(\beta \Lambda)=\beta^L \nu(\Lambda)$. For convenience we will write
$\nu$ for $\nu(\Lambda)$.
Upon rewriting (\ref{eqn-rate-0}) for the scaled lattice, we obtain
\begin{eqnarray}
R_0 & = & h(p)-(1/L) \log_2(\nu(\beta\Lambda)) \nonumber \\
 & = & h(p)-(1/L)\log_2(\beta^L \nu).
\label{eqn-asymp-1}
\end{eqnarray}
Similarly, an expression for $R$ is obtained by rewriting (\ref{eqn-rate-1})  in
order to get
\begin{eqnarray}
R & = & h(p)-(1/L) \log_2(N\nu(\beta\Lambda)) \nonumber \\
 & = & h(p)-(1/L)\log_2(N\beta^L \nu).
\label{eqn-asymp-2}
\end{eqnarray}
From (\ref{eqn-d0}), the two-channel distortion with the scaled lattice $\beta\Lambda$ is
given by
\begin{equation}
\bar{d}_0=G(\Lambda)\nu^{2/L}\beta^{2},
\label{eqn-asymp-3}
\end{equation}
where we have used the fact that $G(\Lambda)=G(\beta\Lambda)$.

In terms of $R_0$, the two-channel distortion $\bar{d}_0$ is thus given by
\begin{equation}
\bar{d}_0=G(\Lambda)2^{2h(p)}2^{-2R_0}.
\label{eqn-asymp-4}
\end{equation}
Now let $N=2^{L(aR+1)}$. Then $R_0=R(1+a)+1$ and
\begin{equation}
\bar{d}_0=\frac{G(\Lambda)2^{2h(p)}}{4} 2^{-2R(1+a)}.
\label{eqn-asymp-5}
\end{equation}
For given $N$ and $R$, the correct scale factor is obtained by solving (\ref{eqn-asymp-3}) for $\beta$:
\begin{equation}
\beta^L=\frac{2^{Lh(p)}2^{-LR(1+a)}}{2^{L}\nu}.
\label{eqn-asymp-6}
\end{equation}
Consider $\tilde{d}$ defined by
\begin{equation}
\tilde{d}=\frac{1}{4N}\sum_{\lambda \in V_0(0)}l^2(e) \beta^2,
\label{eqn-asymp-7}
\end{equation} 
which is the common term in the bounds for $\bar{d}_s$ given in (\ref{eqn-a2-5}). 
The quantity $\beta^2$ arises
because we use the scaled lattices $\beta\Lambda$ and $\beta\Lambda'$. 
It is understood that $e$ is the edge that labels $\lambda$.
The edges in question  in (\ref{eqn-asymp-7}) (the edges in $\calE_d(0)$) have been obtained by choosing the $N$
shortest edges in $\beta\Lambda'$, with one endpoint at $0$ and then shifting these edges so that the
midpoint is as close to the origin as possible. Thus each $l^2(e)$ is of the form $iN^{2/L}/L$ for some
$i\in \ZZ$. 
The term $N^{2/L}$ is a scale factor that comes from the fact that $[\Lambda:\Lambda']=N$ and we
normalize by $L$ because we are working with normalized square lengths.
Let the largest value of $l^2(e)$ in (\ref{eqn-asymp-7}) be equal to $KN^{2/L}/L$ and let $B_i$ be the
number of $l^2(e)$'s that are equal to $iN^{2/L}/L$. 
Then
\begin{equation}
\tilde{d}=\frac{1}{4N}\sum_{i=0}^{K}\frac{\beta^2iB_iN^{2/L}}{L}.
\label{eqn-dom1}
\end{equation}
Our construction of the set $\calE_d(0)$ implies that the $B_i$ can be obtained in terms of the
coefficients of the theta series of the lattice $\Lambda$. To be specific, if
$\Theta_{\Lambda}(z)=\sum_{i}A_iz^{i}$ is the theta series\footnote{$A_i$ is the number of
lattice points $\lambda$ with $L\|\lambda\|^2=i$.} for the lattice $\Lambda$, then we can assert
that
\begin{eqnarray}
B_i & = & A_i,~~~ 0 \leq i <K ~, \nonumber \\
B_K & \leq & A_K.
\end{eqnarray}
This fact will be used a little later.

Now substitute $\beta$ from (\ref{eqn-asymp-6}) and use the fact that $2^{-2aR}=4/N^{2/L}$,
in order to obtain
\begin{eqnarray}
\tilde{d} & = & \frac{2^{2h(p)}2^{-2R(1-a)}}{\nu^{2/L} N^{(1+2/L)} L} \sum_{i=0}^{K}iB_i \nonumber
\\
& = & \Xi \sum_{i=0}^KiB_i,
\label{eqn-dom2}
\end{eqnarray}
where we have defined $\Xi=2^{2h(p)}2^{-2R(1-a)}/(\nu^{2/L} N^{(1+2/L)} L)$. 
The term $\sum_{i=0}^KiB_i$ can be bounded in terms of the $A_i$'s by
\begin{equation}
\sum_{i=0}^{K-1}iA_i \leq \sum_{i=0}^K iB_i \leq \sum_{i=0}^{K}iA_i. 
\end{equation}
Upon defining $S(m)=\sum_{i=0}^m A_i$ and using Abel's summation formula we obtain
\begin{equation}
\sum_{i=0}^{m}iA_i=mS(m)-\sum_{n=0}^{m-1}S(n).
\end{equation}
The term $S(n)$ is the number of lattice points in the first $n$ shells of the lattice
$\Lambda$.  This is roughly the ratio of the volume of $S_L(n^{1/2})$, a sphere of radius
$n^{1/2}$, to $\nu$, the volume of the Voronoi cell of $\Lambda$. To be specific, if
$B_L$ denotes the volume of a sphere of unit radius in $\Reals^L$ then $\lim_{n \rightarrow
\infty} S(n)/n^{L/2}=B_L/\nu$. Thus $S(n)$ is given by
\begin{equation}
S(n)=\frac{B_Ln^{L/2}}{\nu}(1+o(1)),
\label{eqn-numpoints}
\end{equation}
where $\lim_{n\rightarrow \infty}o(1)=0$.
Upon using (\ref{eqn-numpoints}) in order to estimate $S(m)$, we obtain
\begin{eqnarray}
\sum_{i=0}^{m}iA_i & = & \left(\frac{mB_Lm^{L/2}}{\nu}
-\sum_{n=0}^{m-1}\frac{B_Ln^{L/2}}{\nu}\right) (1+o(1))\nonumber \\
& = & \frac{B_L}{\nu}\left(m^{L/2+1}-\sum_{n=0}^{m-1}n^{L/2} \right)(1+o(1))\nonumber \\
& = &
\frac{B_L}{\nu}\left(m^{(L/2+1)}-\frac{m^{(L/2+1)}}{(L/2+1)}+o(m^{L/2+1})\right)(1+o(1))\nonumber \\
& = & \frac{B_L}{\nu}\frac{L}{L+2}m^{(L/2+1)} (1+o(1))
\label{eqn-fromabel}
\end{eqnarray}
On substituting (\ref{eqn-fromabel})  into (\ref{eqn-dom2}) and observing from
(\ref{eqn-numpoints}) that $N=(B_LK^{L/2}/\nu)(1+o(1))$ we obtain
\begin{equation}
\tilde{d}=\frac{1}{B_L^{2/L}(L+2)}2^{2h(p)}2^{-2R(1-a)}(1+  o(1)).
\label{eqn-almost-there}
\end{equation}
But 
\begin{equation}
B_L=\frac{\pi^{L/2}}{\Gamma(L/2+1)},
\end{equation}
and
$G(S_L)$, the normalized second moment of a sphere
in $L$ dimensions is given by
\begin{equation}
G(S_L)=\frac{1}{(L+2)\pi}\Gamma(L/2+1)^{2/L}.
\label{sphere-2m}
\end{equation}
Thus $B_L$ is given in terms of $G(S_L)$ by
\begin{equation}
B_L^{2/L}=\frac{1}{G(S_L)(L+2)},
\end{equation}
and from (\ref{eqn-almost-there}) it follows that
\begin{equation}
\tilde{d}=G(S_L)2^{2h(p)}2^{-2R(1-a)}(1+o(1)).
\end{equation}
The other terms in (\ref{eqn-a2-5}) are $\bar{d}_0$ and $r^*(\Lambda')^2$. The term $\bar{d}_0$
decays as $2^{-2R(1+a)}$ and $r^*(\Lambda)^2$ decays like $\beta^2 N^{2/L}$, which
in turn decays as $2^{-2R}$. Thus $\lim_{R\rightarrow \infty} \bar{d}_s
2^{2R(1-a)}=\lim_{R\rightarrow \infty} \tilde{d}2^{2R(1-a)}$ and we have obtained
our final result:
\begin{equation}
\lim_{R\rightarrow \infty}\bar{d}_s 2^{2R(1-a)}=G(S_L)2^{2h(p)}.
\label{eqn-final-side}
\end{equation}

We end with a comparison with the multiple description rate
distortion bound (\ref{eqn-ecmdsq1}), by letting $L\rightarrow
\infty$ in  (\ref{eqn-asymp-5}) and (\ref{eqn-final-side}).
It is believed that as $L\rightarrow \infty$,
$G(\Lambda)\rightarrow 1/2\pi e$ and it is easily shown from
(\ref{sphere-2m}) that $\lim_{L\rightarrow \infty}G(S_L) = 1/2\pi
e$. Thus our
constructions  are optimal.

\section{Summary and Conclusions}
\label{sec-summary}
The problem of lattice vector quantizer design is addressed for the two-channel multiple
description. The main problem in the design, a labeling problem, is solved.
A systematic construction technique is developed which is suitable for general lattices. 
Specific constructions have been provided for $A_2$ and $\ZZ^i$, $i=1,2,4,8$. 
Finally an asymptotic analysis reveals that
performance arbitrarily close to the multiple description rate distortion bound can be
obtained.

Open issues related to this work are detailed constructions for other
lattices, extensions to the asymmetric case and extensions to greater than two descriptions.

\vspace{0.3in}
\begin{center}
{\bf Acknowledgement}\\
\end{center}
We thank Gabriele Nebe for suggesting the group of order 16 that we used for the
lattice $\ZZ^8$. We also thank the reviewers for their detailed and thoughtful
comments and suggestions.

\clearpage

\end{document}